\PassOptionsToPackage{sort&compress}{natbib}
\PassOptionsToPackage{dvipsnames}{xcolor}
\documentclass[a4paper, preprint, 3p]{elsarticle}

\usepackage[T1]{fontenc}
\usepackage[utf8]{inputenc}
\usepackage[USenglish]{babel}
\usepackage{derivative}
\usepackage{sansmath}

\usepackage{amssymb}
\usepackage{amsthm}
\usepackage{thmtools}
\usepackage{mathtools}
\usepackage{bm}
\usepackage{algorithm}
\usepackage[noEnd=false]{algpseudocodex}
\usepackage{graphicx}
\usepackage{float}
\usepackage{natbib}
\usepackage{amsmath}
\usepackage{wrapfig}
\usepackage{caption}
\usepackage{subcaption}
\usepackage{graphicx}
\usepackage{multirow}
\usepackage{csvsimple}
\usepackage[plainpages=false,pdfpagelabels,hidelinks,unicode]{hyperref}
\usepackage{soul} 
\usepackage{booktabs}
\usepackage{csquotes}
\usepackage{xspace}
\usepackage{adjustbox}
\usepackage{xcolor}
\declaretheoremstyle[
  bodyfont=\normalfont\itshape,
  headformat=\NAME\ \NUMBER\NOTE,
]{myplain}
\declaretheoremstyle[
  headformat=\NAME\ \NUMBER\NOTE,
]{mydefinition}
\newcommand{\envqed}{{\lower-0.3ex\hbox{$\triangleleft$}}}
\declaretheorem[style=myplain,numberwithin=section]{theorem}

\declaretheorem[style=mydefinition,numberlike=theorem,qed=\envqed]{remark}

\usepackage[frozencache,cachedir=minted_cache-secure_simulations__arXiv]{minted}
\setminted{
    autogobble,
    fontsize=\small,
    frame=none,
    numbers=left
}

\DeclareMathOperator{\circshift}{circshift}
\DeclareMathOperator{\rotate}{rotate}

\usepackage{todonotes}
\setuptodonotes{inline, caption={}}
\algrenewcommand\textproc{}

\newcommand{\orcid}[1]{ORCID:~\href{https://orcid.org/#1}{#1}}

\graphicspath{{..}{figures/}{talk/figures/}}

\date{\today}

\makeatletter
\hypersetup{pdfauthor={Arseniy Kholod and Yuriy Polyakov and Michael Schlottke-Lakemper}}
\hypersetup{pdftitle={Secure numerical simulations using fully homomorphic encryption}}
\makeatother

\begin{document}

\begin{frontmatter}

\title{Secure numerical simulations using fully homomorphic encryption}

\author[1,2,3]{Arseniy Kholod\fnref{orcidAK}}
\fntext[orcidAK]{\orcid{0009-0002-9457-4867}}

\author[4]{Yuriy Polyakov\fnref{orcidYP}}
\fntext[orcidYP]{\orcid{0000-0002-5566-3763}}

\author[1,2,3]{Michael Schlottke-Lakemper\fnref{orcidMSL}}
\fntext[orcidMSL]{\orcid{0000-0002-3195-2536}}
\ead{michael.schlottke-lakemper@uni-a.de}
\cortext[cor1]{Corresponding author}

\address[1]{High-Performance Scientific Computing, Institute of Mathematics, University of
Augsburg, Germany}

\address[2]{Centre for Advanced Analytics and Predictive Sciences, University of Augsburg, Germany}

\address[3]{Applied and Computational Mathematics, RWTH Aachen University, Germany}

\address[4]{Duality Technologies, Inc., United States}

\begin{abstract}
    Data privacy is a significant concern when using numerical simulations for sensitive
information such as medical, financial, or engineering data---especially in untrusted
environments like public cloud infrastructures. Fully homomorphic encryption (FHE) offers a
promising solution for achieving data privacy by enabling secure computations directly on
encrypted data. Aimed at computational scientists, this work explores the viability of
FHE-based, privacy-preserving numerical simulations of partial differential equations. The
presented approach utilizes the Cheon-Kim-Kim-Song (CKKS) scheme, a widely used FHE method for approximate
arithmetic on real numbers. Two Julia packages are introduced, OpenFHE.jl and
SecureArithmetic.jl, which wrap the OpenFHE C++ library to provide a convenient interface
for secure arithmetic operations. With these tools, the accuracy and performance of key FHE
operations in OpenFHE are evaluated, and implementations of finite difference
schemes for solving the linear advection equation with encrypted data are demonstrated. The
results show that cryptographically secure numerical simulations are possible, but that
careful consideration must be given to the computational overhead and the numerical errors
introduced by using FHE. An analysis of the algorithmic restrictions imposed by FHE
highlights potential challenges and solutions for extending the approach to other models and
methods. While it remains uncertain how broadly the approach can be generalized to more complex
algorithms due to CKKS limitations, these findings lay the groundwork for further research on
privacy-preserving scientific computing.

\end{abstract}

\begin{keyword}
  secure numerical simulations \sep
  fully homomorphic encryption \sep
  CKKS scheme \sep
  OpenFHE \sep
  Julia programming language

  \MSC[2020]
  65-04 \sep 
  65M06 \sep 
  65Y99 
\end{keyword}

\end{frontmatter}

\section{Introduction}\label{sec:introduction}

Partial differential equations (PDEs) are used to model phenomena across
scientific fields ranging from physics and engineering to biology and finance.
Since many PDEs
cannot be solved analytically, numerical methods are used to approximate their
solutions in scientific, industry, and business applications.
In certain cases, the input data and simulation
results involve highly sensitive information—--such as personal health records, financial
details, or proprietary engineering designs--—that must be appropriately safeguarded. For
instance, finite element simulations of the aorta used in personalized approaches for
treating cardiovascular disease rely on patient-specific data \cite{djukic2019numerical,
perrin2015patient}. These datasets must be handled securely to comply with regulations like the
General Data Protection Regulation (GDPR) in the European Union or the Health Insurance
Portability and Accountability Act (HIPAA) in the United States, posing challenges for
scientific research as well as practical implementations \cite{clause2004conforming,
conley2018gdpr}.
Another example arises in population genetics, where the evolution of gene frequencies
has long been modeled by reaction-diffusion equations
\cite{kimura1964diffusion, vlad2004enhanced, lou2013pdemodels} that are solved numerically \cite{gutenkunst2009inferring, andreianov2011analysis,
zhao2013complete}.
A practical issue here is working with combined datasets from multiple providers (which do not trust each other), while
protecting the privacy of individual genetic information \cite{gymrek2013identifying,
brenner2013bigleak, BGP20}.

To safeguard data, encryption is commonly employed both at rest and during transmission.
However, once the data is decrypted for processing, it becomes vulnerable to attacks.
The data privacy risks associated with decryption for processing become
pronounced when data must be processed remotely, as is often the
case with public cloud computing platforms. Limited local resources or the need for
distributed systems may necessitate offloading computations to external cloud infrastructures,
thereby exposing decrypted data to potential attacks. With the rapid expansion of cloud
computing \cite{GansHerveEtAl23}, ensuring data privacy during computation has emerged as a
pressing concern for both practitioners and researchers \cite{acar2018survey, MarcollaSucasasEtAl22}.

A possible solution to this conundrum is the use of what is known as \emph{fully homomorphic
encryption} (FHE).
With FHE, it is possible to perform computations directly on the
encrypted data.
That is, the result of an operation will be same as if the data had been
decrypted first, processed, and then encrypted again. Since in case of FHE the data is never
available in plaintext, even a fully compromised compute system would not leak any sensitive
information. The notion of homomorphic encryption (homomorphisms) was introduced by Rivest, Adleman, and
Dertouzos, who in 1978 suggested that such a scheme might be possible
\cite{Rivest1978}. The first viable scheme was proposed much later by Gentry in 2009
\cite{Gentry09},
when he showed that fully homomorphic encryption is achievable using
lattice-based cryptography. Since then, several FHE schemes have been proposed, with the
Cheon-Kim-Kim-Song (CKKS) scheme being the most popular option for computations with real
numbers \cite{CheonKimEtAl17}.
Other commonly used schemes include the Brakerski/Fan--Vercauteren
(BFV)~\cite{Brakerski12, FanVercauteren12} and Brakerski--Gentry--Vaikuntanathan
(BGV)~\cite{BrakerskiGentryEtAl14} schemes for computations with integers, and the Ducas--Micciancio (DM/FHEW)~\cite{DM15} and
Chillotti--Gama--Georgieva--Izabachene (CGGI/TFHE)~\cite{ChillottiGamaEtAl19} schemes for evaluating Boolean and small-integer arithmetic circuits. A good overview of FHE schemes and their applications can be
found in \cite{MarcollaSucasasEtAl22}.

Fully homomorphic encryption schemes represent an active field of research, and
applications of FHE schemes are used in disciplines such as privacy-preserving
machine learning \cite{LeeKangEtAl22}, analysis of encrypted medical or genomic data
\cite{BosLauterEtAl14, ChenGilad-BachrachEtAl18, BlattGusevEtAl20, BGP20}, or for processing data
in financial services \cite{MastersHuntEtAl19}.
One area where FHE has yet to be applied
is the field of computational physics. There are mainly three reasons for
this: First, the concept of FHE is still new and not well known outside of the
cryptography community. Second, while there are some excellent open-source libraries
available for FHE, they currently require a high level of expertise to use.
Third, since FHE only supports a very limited subset of the usual arithmetic operations,
implementing secure numerical algorithms is not straightforward and requires a good
understanding of FHE.
In this paper, we aim to address these issues
by examining fully homomorphic encryption from a scientific computing perspective
and introducing it to the computational physics community.

Two representative scenarios highlight the practical relevance of FHE-secured numerical
simulations.
Both involve a data
owner (holding sensitive model parameters or input data) and a compute
provider (offering computational resources and/or proprietary simulation software).
In the first scenario, the data owner must ensure the confidentiality of the input data, for
example, a medical researcher running patient-specific simulations on a public cloud
platform.
In the second scenario, both parties require confidentiality: the input data must remain
private, and the compute provider must protect intellectual property in the simulation code.
One example is the use of proprietary multiphysics solvers to analyze sensitive design data in aerospace or energy applications.
In both scenarios, the workflow proceeds as follows: the data owner locally encrypts the
input using a homomorphic encryption scheme such as CKKS and transmits the ciphertexts to
the compute provider. The provider runs the numerical simulation directly on the encrypted
data and returns the encrypted results. Decryption and post-processing are
then performed locally by the data owner.

To explore the feasibility of secure simulations in a scientific computing context, this paper is structured as follows:
Sec.~\ref{sec:fhe} begins with a brief overview of the CKKS scheme for computational
scientists. Moreover, we introduce the open-source C++ library OpenFHE and the Julia
packages OpenFHE.jl and SecureArithmetic.jl, which provide a convenient interface for secure
arithmetic operations in Julia for rapid experimentation. We then use these packages in
Sec.~\ref{sec:accuracy_performance} to analyze the performance and accuracy of basic FHE
operations to establish a baseline for more complex algorithms. In
Sec.~\ref{sec:secure_simulation}, we implement several finite difference schemes to solve
the linear advection equation using FHE. They serve as prototypes for FHE-secured numerical
simulations, and their performance and accuracy are analyzed in Sec.~\ref{sec:results}. While secure numerical simulations through FHE are feasible, the
computational overhead is currently significant for many practical scenarios, and the
encryption-induced numerical errors must be taken into account when designing applications.
Sec.~\ref{sec:extension} explores challenges and potential solutions for extending secure
numerical simulations to other mathematical models and numerical methods. Finally, we
summarize our findings in Sec.~\ref{sec:conclusions} and provide an outlook on future work.
All code used in this paper, as well as the data generated during the experiments, is
available in our reproducibility repository \cite{reproKholodSchlottkeLakemper24}.

\section{Fully homomorphic encryption}\label{sec:fhe}

\begin{figure}[!htb]
    \centering
    \includegraphics[width=0.85\textwidth]{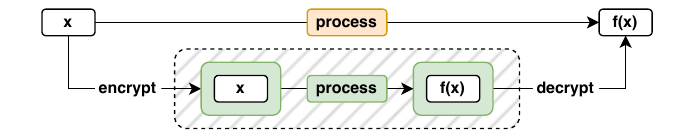}
    \caption{Comparison of data processing in the clear (top path) to processing it under homomorphic encryption (bottom path). Everything within the hatched area is secure, allowing one to process sensitive information in untrusted environments.}
    \label{fig:he-overview}
\end{figure}

The fundamental concept of homomorphic encryption is to perform operations on encrypted data
without decrypting it first.
This is achieved by using a special encryption scheme that allows for the evaluation of
arbitrary functions on encrypted data, such that the result---once decrypted---is the same
as if the data had been processed in the clear.
Thus, in a sense, a homomorphic encryption scheme is transparent to the
application of arbitrary functions to the data. This is illustrated in
Fig.~\ref{fig:he-overview}, where some data $x$ is first homomorphically encrypted, then
processed into $f(x)$, and finally decrypted (bottom path). This yields the same result as
if the data had been directly processed in the clear (top path).

With \emph{partially homomorphic encryption}, only one type of arithmetic operation can be
performed on the encrypted data, i.e., either addition or multiplication.
\emph{Leveled homomorphic encryption} allows both addition and multiplication, supporting
arbitrary computation circuits of bounded (pre-determined) depth.
Finally, \emph{fully
homomorphic encryption} (FHE) enables the evaluation of arbitrary functions of unbounded
complexity on encrypted data, and is thus the most powerful type of homomorphic encryption.

We begin with a high-level overview of the CKKS scheme for fully homomorphic encryption of
real numbers, followed by a brief introduction to the software packages developed for and
used in this work. The discussion is intended for computational scientists unfamiliar with
cryptography. A more comprehensive treatment of the underlying principles can be found
in~\cite{Gentry10}.

\subsection{Secure approximate arithmetic with the CKKS scheme}\label{sec:ckks}

The CKKS (Cheon-Kim-Kim-Song) scheme \cite{CheonKimEtAl17} is a fully homomorphic encryption
method for approximate arithmetic on real numbers, built upon the BGV scheme for FHE with
integers \cite{BrakerskiGentryEtAl14}. It obtains its security from the
\emph{ring learning with errors} (RLWE) problem, which is a specialization of the
\emph{learning with errors} problem \cite{Regev09} for polynomial rings over finite fields
\cite{LyubashevskyPeikertEtAl10}.
Like other FHE methods,
CKKS can be used either as a symmetric encryption scheme (a single key for encryption and decryption) or an
asymmetric encryption scheme
(a public key for encryption and a separate private key for
decryption). In this manuscript, we solely focus on its use in the asymmetric
encryption context.

\begin{figure}[!htb]
    \centering
    \includegraphics[width=0.85\textwidth]{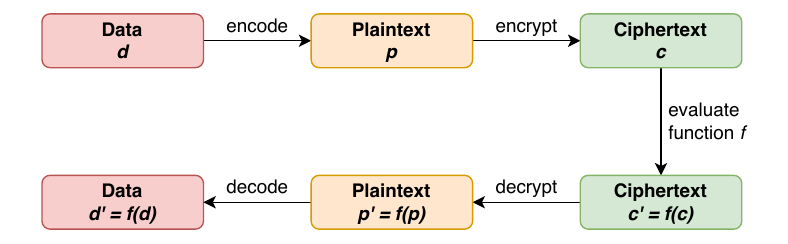}
    \caption{CKKS overview (following \cite{PanChaoEtAl24}).}
    \label{fig:fhe-overview}
\end{figure}

Fig.~\ref{fig:fhe-overview} outlines the general procedure when using CKKS for FHE. We start
with the raw user data $d$, which in case of the CKKS scheme is typically a vector of
double-precision floating point numbers (complex numbers are also supported by the CKKS scheme, but they are rarely used in practice). This data is first mapped to a suitable polynomial
representation, called the \emph{plaintext} $p$. The plaintext can then be encrypted using
the public key to obtain the \emph{ciphertext} $c$. The ciphertext $c$ is now ready to be
processed by some arbitrary function $f$, which yields another ciphertext $c' = f(c)$.
Finally, the result can be decrypted using the private key to obtain the plaintext
representation $p' = f(p)$, and then decoded back to a vector of real numbers $d' = f(d)$,
which is the desired output.

The encrypted ciphertext is obtained by adding a random mask generated via RLWE to the plaintext in a
controlled manner, effectively hiding the raw data behind a veil of randomness. The computational security of the CKKS scheme, just like any other FHE scheme based on RLWE, relies on the
mathematical complexity of solving a noisy high-dimension linear system of equations. This task remains
computationally infeasible, even for quantum computers. However, in CKKS the noise also introduces a
small error in the computations that is not fully recoverable during the decryption
procedure. In this sense, the CKKS scheme is an approximate arithmetic method (which distinguishes it from exact FHE schemes, such as BGV and BFV), and the results of
CKKS computations are always subject to some level of approximation error. The magnitude of this error
depends on the parameters of CKKS and the number and type of operations performed while
evaluating the function $f$.

In terms of arithmetic operations that can be evaluated on the ciphertext, the CKKS scheme
supports element-wise \emph{addition} and \emph{multiplication} with ciphertexts,
plaintexts, or scalar values, as well as \emph{rotation} operations.
The rotation operation performs a circular shift of the data elements in the ciphertext.
That is, it moves the data elements by a fixed number of positions while preserving their
order, with elements that are shifted out of one end being reintroduced at the opposite end.
Rotations are essential for implementing higher-level operations such as vector summation or matrix-vector products.
From these three basic operations, more
complex functions can be constructed.
Each of these operations introduces a small amount of noise into the ciphertext, which
slightly degrades the accuracy of the underlying plaintext data.
Moreover, multiplication substantially increases the internal noise level of the ciphertext, to the point where it can eventually prevent correct decryption. To avoid this, a process known as \emph{rescaling} is applied to keep the noise level in check.
Due to this, and because of how the CKKS scheme
is constructed, the multiplicative depth of computations in CKKS, i.e., the number of consecutive multiplications, is limited to a fixed
number of levels. To overcome this limitation, a technique called \emph{bootstrapping} can
be applied after all levels of a ciphertext are used. This bootstrapping procedure resets the
available multiplicative depth
in the ciphertext, enabling more computations. Since bootstrapping
may be used repeatedly, it effectively removes the restrictions on the algorithmic depth of
computations in CKKS (for many practical applications where approximation error growth is limited).

\begin{remark}
    With the availability of addition, multiplication, and rotation as arithmetic
    operations, the CKKS scheme is indeed capable of evaluating arbitrary functions on
    encrypted data. However, the restriction to just these three basic operations means that, in
    practice, iterative procedures or polynomial approximations are required to evaluate even moderately complex
    functions. For example, to determine the multiplicative inverse of a number (e.g., for a
    division operation), variants of the Newton-Raphson method \cite{CetinDorozEtAl15} or
    Goldschmidt's algorithm \cite{MoonOmarovEtAl24} are often used.
\end{remark}

In the following sections, we discuss the CKKS scheme and its main operations in more
detail. For clarity, scalar quantities (including polynomials) are written in regular font,
while vectors are denoted by bold font. While not strictly necessary to understand the
remainder of the paper, these sections provide a useful background for the subsequent
sections, especially for readers unfamiliar with FHE methods.
Good high-level descriptions
of the CKKS scheme can also be found in, for example, \cite{Huynh20,Polyakov21, Inferati22}.

\subsubsection{Plaintext encoding and decoding}\label{sec:encoding}

The encoding mechanism converts the user data from a representation as a vector of
double-precision floating point numbers into a polynomial representation. Polynomials in the
CKKS scheme are \emph{cyclotomic polynomials} $\mathcal{R} = \mathbb{Z}[X]/(X^{N_R}+1)$ with
$N_R$ integer coefficients, where $N_R$ is always a power of two and is called the
\emph{ring dimension}. The \emph{batch size} or \emph{capacity}, i.e., the number of slots
available for packing data into a plaintext, is also a power of two and can be at most as
large as half the ring dimension $N_R/2$.

For encoding, the input vector is treated as coefficients of a polynomial. This polynomial is first converted to a new representation by evaluating it at complex roots of
unity of the $N_R$-th cyclotomic polynomial using a procedure
resembling
an inverse Fast
Fourier Transform (FFT).
This enables component-wise vector multiplication using FHE. The result of the inverse FFT
is then scaled by a large value, the \emph{scaling factor} (or scaling modulus),
and rounded to the nearest integer values to make it compatible with RLWE.
At this point, the
data is still not not encrypted. Due to the limited precision of (double) floating-point numbers and rounding
error, the encoding step already introduces a small error into the data.

The decoding mechanism downscales the integers back to floating-point values and then 
converts the polynomial back into a data vector by evaluating the polynomial at the complex roots of unity,
an operation similar to forward FFT. The noise increase brought about by this operation is
typically negligible as the error introduced by the decoding operation is usually smaller
than the existing approximation error resulting from prior CKKS operations, such as
encryption and computations.

\subsubsection{Levels and ciphertext modulus}\label{sec:levels}

CKKS is a leveled FHE scheme (with bootstrapping), where a freshly encrypted ciphertext
starts with a level $l > 0$. With each multiplication of the ciphertext, the
level $l$ is decreased by one until $l = 0$. Once the level reaches zero, no more
multiplications are possible with this ciphertext. Bootstrapping, as discussed in
Sec.~\ref{sec:bootstrapping}, alleviates this restriction by resetting the level to a value
greater than zero. Addition or rotation operations typically do not consume levels.

Therefore, we usually set FHE parameters using the desired multiplicative depth of the
computation. We denote the maximum multiplicative depth as $l_\text{max}$. Based on the
multiplicative depth, precision requirements, and security level, two main lattice
parameters are configured: ring dimension $N_R$ (as described in Section~\ref{sec:encoding})
and \emph{ciphertext modulus} $Q$.

The ciphertext modulus gets reduced
after each multiplication by the \emph{rescaling operation}.
This operation discards the least significant part of the ciphertext---analogous to
truncation in floating-point arithmetic---thereby controlling noise growth and maintaining a
fixed scale for subsequent operations.
Its purpose
is explained later in Sec.~\ref{sec:multiplication}. It is thus convenient to use the
notation of $Q_l$, where $l$ is the current level. We start with $l = l_\text{max}$, i.e.,
$Q=Q_{l_\text{max}}$.

At a more granular level, the ciphertext modulus is a product of $l + 1$ factors and is
represented as $Q_l = q_0 \cdot q_1 \cdot q_2 \dots q_{l-1} \cdot q_{l}$. The rescaling
operation reduces
the ciphertext modulus from $Q_l$ to $Q_{l-1}$. $q_0$ is the
decryption modulus (often called the \emph{first modulus}), and for $i = 1 \dots l$,
$q_i$ is equal (exactly or approximately) to the scaling factor.

For efficiency, the OpenFHE implementation of CKKS (which we use in this paper, see
also Sec.~\ref{sec:software}) uses the Residue Number System (RNS) representation of large
numbers, with all values $q_i$ being co-prime. For the RNS implementation and a fixed ring
dimension $N_R$, most of FHE operations scale approximately linearly with the level $l$
\cite{OpenFHE}.

\subsubsection{Encryption and decryption}\label{sec:encryption}
CKKS is a public key encryption scheme with two keys, a public and a private one. The public
key $\bm{\mathrm{pk}}$ is used for encryption and can be safely shared, while decryption
uses the private key $\bm{\mathrm{sk}}$. The public key is generated using three
polynomials: a uniformly random polynomial $a$, a small secret polynomial $s$, and a small
error polynomial $e_\mathrm{pk}$. The secret key is defined as
\begin{equation}
    \bm{\mathrm{sk}} = (1, s)
\end{equation}
and the public key is defined as
\begin{equation}
    \bm{\mathrm{pk}} \coloneqq (\mathrm{pk}_0, \mathrm{pk}_1) = (-a \cdot s + e_\mathrm{pk}, a) \pmod{Q_l},
\end{equation}
where $\cdot$ denotes polynomial multiplication. The term $\mathrm{pk}_0 \coloneqq -a \cdot
s + e_\mathrm{pk}$ corresponds to the RLWE problem, which can be informally stated as follows: if $a$
and $\mathrm{pk}_0$ are given, it is computationally hard to find $s$ (and $e_\mathrm{pk}$).
Furthermore, we use $\pmod{Q_l}$ to emphasize where arithmetic operations are
performed modulo the current ciphertext modulus $Q_l$.
The error polynomial $e_\mathrm{pk}$ is essential for
the security of the scheme, as it prevents an attacker from recovering the secret $s$ from
the public key $\bm{\mathrm{pk}}$ by solving the corresponding linear systems of equations.

To encrypt the data previously encoded in the plaintext message $\mu$, which itself is a
polynomial, three small random polynomials are generated: $u$, $e_1$, and $e_2$. The public
key $\bm{\mathrm{pk}}$ is then used to encrypt the encoded data into a ciphertext $\bm{c} =
(c_0, c_1)$ as follows:
\begin{equation}
    \begin{aligned}
        c_0 &= \mathrm{pk}_0 \cdot u + e_1 + \mu = (-a \cdot s + e_\mathrm{pk}) \cdot u + e_1 + \mu \pmod{Q_l}, \\
        c_1 &= \mathrm{pk}_1 \cdot u + e_2 = a \cdot u + e_2 \pmod{Q_l}.
    \end{aligned}
    \label{eqn:encryption}
\end{equation}
As can be seen from Eqn.~\ref{eqn:encryption}, a ciphertext
always consists of a pair of polynomials. Similar to the public key before, the plaintext
message is masked by the pseudorandomness coming from the RLWE problem.

To decrypt the encrypted data $\bm{c} = (c_0, c_1)$, the private key $\bm{\mathrm{sk}}$ is
used to obtain an approximate value of the plaintext message $\tilde{\mu}$ by computing the
inner product $\langle \cdot, \cdot\rangle$ of the ciphertext $\bm{c}$ with $\bm{\mathrm{sk}}$ as
\begin{equation}
    \begin{aligned}
    \tilde{\mu} &= \langle \bm{c}, \bm{\mathrm{sk}}\rangle = c_0 + c_1 \cdot s \pmod{Q_l}\\
                &= \mu + e_\mathrm{pk} \cdot u + e_1 + e_2 \cdot s = \mu + v \approx \mu,
    \label{eqn:decryption}
    \end{aligned}
\end{equation}
where $v = e_\mathrm{pk} \cdot u + e_1 + e_2 \cdot s$ is small compared to the plaintext message $\mu$
by construction. Thus, the decryption operation is not exact and adds a small error to the
decrypted data.

\subsubsection{Addition}\label{sec:addition}
Addition of two ciphertexts is done by adding the corresponding polynomials element-wise,
i.e.,
\begin{equation}
    \bm{c}' = \bm{c} + \hat{\bm{c}} = (c_0 + \hat{c}_0, c_1 + \hat{c}_1) \pmod{Q_l}.
    \label{eqn:addition}
\end{equation}
This operation results in an error growth from $e$ to $e+\hat{e}$ in decrypted data, which
is still small compared to the data itself \cite{CheonKimEtAl17}. To add a plaintext to a ciphertext, the
plaintext is encoded and then extended to the ciphertext space (no additional error is introduced by this operation). Addition of ciphertexts and scalars is done similarly (but no encoding is needed in this case).

\subsubsection{Key Switching}
\label{sec:keyswitching}

Unlike addition, the two other basic arithmetic operations, multiplication and rotation, require executing a special maintenance procedure called \textit{key switching}. This key switching procedure needs an additional public key, which is often called an evaluation key or key-switching hint/key. We use the term evaluation key in the paper.

Key switching is necessary because a multiplication or rotation transforms both the encrypted
message and the underlying secret key. To transform the resulting ciphertext back to the
original secret key, we perform the key switching procedure using a properly generated
evaluation key. The evaluation key is essentially an encryption of the transformed secret
key under the original secret key, which is done in a way that minimizes the noise increase
associated with key switching.

In CKKS, the key switching procedure is more computationally expensive than the
actual multiplication or rotation.
This is because key switching involves multiple number-theoretic transforms (a specialized
version of the discrete Fourier transform), which dominate the overall cost of these
encrypted operations in practice~\cite{BlattGusevEtAl20, yang2022bandwidth}. As such, key switching is
typically the main performance bottleneck in CKKS-based FHE applications.

\subsubsection{Multiplication}\label{sec:multiplication}
Multiplication of two CKKS
ciphertexts $\bm{c}$ and $\hat{\bm{c}}$ is a considerably more complex operation than
addition, which is why we only provide a brief overview here. It requires three steps:
\begin{enumerate}
    \item Multiplication of the polynomial pairs $\bm{c} = (c_0, c_1)$ and $\hat{\bm{c}} =
    (\hat{c}_0, \hat{c}_1)$ to obtain a polynomial triple $\bm{d} = (d_0, d_1, d_2) = (c_0
    \cdot \hat{c}_0, c_0 \cdot \hat{c}_1 + \hat{c}_0 \cdot c_1, c_1 \cdot \hat{c}_1)$. This
    is typically called a \textit{tensor product}.
    \item \emph{Relinearization} of the polynomial triple $d$ to reduce it to polynomial pair
    $\tilde{\bm{c}} = (\tilde{c}_0, \tilde{c}_1)$ again. This step requires a key switching operation.
    \item \emph{Rescaling} of the ciphertext $\tilde{\bm{c}}$ to reduce the message to the same
    scale as before and to reduce the ciphertext modulus from $Q_l$ to $Q_{l-1}$.
\end{enumerate}

As we saw in Sec.~\ref{sec:encryption}, after encryption the first polynomial of the
ciphertext $c_0$ is linear with respect to the secret key $\bm{\mathrm{sk}}$. 
After multiplying two
ciphertexts using the tensor product, we get a quadratic polynomial in secret key polynomial
$s$, i.e., the decryption would be evaluated as $\langle \bm{d}, (1, s, s^2)\rangle = d_0 + d_1
\cdot s +
d_2 \cdot s^2$.
Each
subsequent multiplication will increase the degree of $s$ even further; one can think of a multiplication of two ciphertexts as a product of two corresponding decryption polynomials. This implies the size
of the ciphertext after a multiplication of two ciphertexts of sizes $i$ and $j$ will grow
to $i+j-1$ polynomials.
To maintain the compact representation of ciphertexts, the relinearization procedure
compresses the ciphertext by reducing it by one polynomial, in this case from three to two
polynomials.

Another issue is that the multiplication of ciphertext leads to error growth because
a small error polynomial is multiplied by encrypted message values. This problem is mitigated by the
rescaling technique, which reduces the scaling factor to the one used during
encryption, hence truncating the least significant error that appeared as the result of
multiplication. Each rescaling reduces the ciphertext modulus $Q_l$ by a factor comparable
to the scaling factor, thereby decreasing the ciphertext level $l$ and thus the number of
multiplications still available.

In leveled FHE, one therefore tries to minimize the required multiplicative depth of an
algorithm by evaluating a chain of multiplications using a binary tree approach. For
example, instead of computing $\bm{c'} = \bm{c}_1 \cdot \bm{c}_2 \cdot \bm{c}_3 \cdot
\bm{c}_4$ with three consecutive multiplications (and thus requiring three levels), one can
first compute $\bm{r}_1 = \bm{c}_1 \cdot \bm{c}_2$ and $\bm{r}_2 = \bm{c}_3 \cdot \bm{c}_4$
using one level, and then combine them to $\bm{c'} = \bm{r}_1 \cdot \bm{r}_2$ using a second
level. In the na\"ive case, to multiply $n$ ciphertexts, the multiplicative depth of $n-1$
is needed. If the binary tree method is used instead, the multiplicative depth is reduced to
approximately $\log_2 (n)$ (see also \cite{Polyakov21}).

\subsubsection{Rotation}\label{sec:rotation}
The multiplication and addition operations are element-wise over the entire ciphertext,
i.e., they are always performed between corresponding elements of the user-provided data
vectors. However, many practical algorithms require interactions between elements within the
same vector or access to specific indices, which goes beyond basic element-wise
operations.
In such cases, the rotation operation is used, which cyclically shifts encrypted data by
some rotation index. When the rotation is evaluated, the underlying secret key is also
implicitly rotated.
To change the underlying secret key back to the original one (to support further homomorphic
computations), one needs to apply the key switching operation. The rotation itself
is cheap (effectively it is just reindexing), but the required key switching operation
following it is computationally expensive.

An important property of
rotation is that it is circular with respect to the batch size/capacity.

\subsubsection{Bootstrapping}\label{sec:bootstrapping}
As described in Sec.~\ref{sec:levels} and \ref{sec:multiplication}, the necessity to rescale
the ciphertext limits the number of multiplicative operations that can be performed, since
no more multiplications are possible once the level $l$ reaches zero. The bootstrapping
operation alleviates this restriction by resetting the level to $l > 0$. It involves
approximating a homomorphic decryption procedure on a ciphertext using an encrypted version
of the secret key (this is done indirectly in CKKS), and then re-encrypting it to obtain a
refreshed ciphertext. In theory, the number of multiplications that can be performed with
a ciphertext is unlimited when using bootstrapping. However, bootstrapping also adds
a significant approximation error to the ciphertext and is computationally very expensive.
Due to the high complexity of the bootstrapping operation, we do not go into more detail here.
In-depth descriptions of various bootstrapping implementations are covered elsewhere~\cite{Bootstrapping}. 

\subsection{Software libraries for fully homomorphic encryption}\label{sec:software}

While the CKKS scheme is a powerful method for FHE with real numbers, it is also an advanced
technique that requires a deep understanding of the underlying
mathematical methods to be used effectively.
It is therefore advisable to use a software
library that provides an implementation of the CKKS scheme and other necessary operations.
Several such libraries are available, e.g., SEAL \cite{SEAL23}, HElib
\cite{HaleviShoup20}, and OpenFHE \cite{OpenFHE}.

In this paper, we focus on the OpenFHE library, which is actively maintained and includes most recent CKKS optimizations~\cite{KimPapadimitriouEtAl20}. We also present two Julia packages OpenFHE.jl
and SecureArithmetic.jl, which provide a convenient interface to OpenFHE in the Julia
programming language. We now give a brief overview of these tools.

\subsubsection{OpenFHE}\label{sec:openfhe}

OpenFHE \cite{OpenFHE} is an open-source FHE library that offers efficient C++
implementations of all common FHE schemes, including the CKKS scheme with bootstrapping. 
OpenFHE contains many examples that are especially
useful for new users when designing their applications. The OpenFHE library puts an emphasis
on usability. For example, the relinearization and rescaling operations after ciphertext
multiplication are handled automatically, allowing the library to be used also by non-experts in FHE.
This can be observed in a short code snippet in Listing~\ref{lst:openfhe-vs-openfhejl}
(left), where the CKKS internals for multiplication are hidden from the user. Some
code for initializing the setup was omitted for clarity, e.g., to create the cryptographic
context object (named \texttt{cc} in this listing) or the public and private keys. For more
details about OpenFHE, please refer to \cite{OpenFHE}.

\subsubsection{OpenFHE.jl}\label{sec:openfhejl}

The Julia programming language \cite{BezansonEdelmanEtAl17} is designed for technical
computing, with a simple, expressive syntax and high computational performance. It combines
the ease of use of high-level languages like Python with the speed of compiled languages
such as C or Fortran, making it ideal for scenarios that require both rapid prototyping and
fast execution. Julia comes with its own package manager, which allows users to easily
install and manage additional packages, including those written in other languages.

It is this feature that we leverage in our OpenFHE.jl package
\cite{schlottkelakemper2024openfhejulia}, which provides a Julia wrapper for the OpenFHE
library. When installing OpenFHE.jl, pre-built binaries of the OpenFHE library are
automatically downloaded without the user having to compile anything locally. The C++
functionality is exposed in Julia via the CxxWrap.jl package \cite{Janssens20}.
Besides offering a native Julia interface to the OpenFHE library, OpenFHE.jl does not
provide any extra functionality. It is intended to be used as a backend for other Julia
packages that require FHE capabilities. As can be seen in
Listing~\ref{lst:openfhe-vs-openfhejl} (right), the syntax of OpenFHE.jl is very similar to
that of OpenFHE, making it easy to port code between the two programming languages.

\subsubsection{SecureArithmetic.jl}\label{sec:securearithmeticjl}

To make the use of FHE more accessible to non-experts in Julia, we created the
SecureArithmetic.jl package \cite{schlottkelakemper2024securearithmetic}.
It inherits the low-level
functionalities of the OpenFHE.jl package, but provides them through a more
convenient, higher-level Julia interface.
Using the SecureArithmetic.jl package, users can write arithmetic
expressions in a common mathematical notation, significantly simplifying the code.
Furthermore, SecureArithmetic.jl supports the ability to use the same code for secure and
non-secure computations, which is especially useful during the design stage. It allows one
to debug FHE algorithms without having to go through the actual encryption and decryption
steps, making the execution much faster during the prototyping phase.
Listing~\ref{lst:securearithmeticjl} shows an example of how the ciphertext multiplication
is performed using SecureArithmetic.jl. Compared to OpenFHE/OpenFHE.jl in
Listing~\ref{lst:openfhe-vs-openfhejl}, SecureArithmetic.jl provides a more high-level
interface that is closer to the mathematical notation. However, since it is built on top of
OpenFHE.jl and OpenFHE, it is still possible to access the lower-level functions if needed.

\begin{listing}[!htbp]
    \centering
    \begin{minipage}[t]{0.45\textwidth}
        \begin{minted}[fontsize=\small, frame=none]{cpp}
            std::vector<double> v = {1.0, 2.0, 3.0, 4.0};
            Plaintext p = cc->MakeCKKSPackedPlaintext(v);
            auto c = cc->Encrypt(public_key, p);
            auto c_squared = cc->EvalMult(c, c);
            Plaintext result;
            cc->Decrypt(private_key, c_squared, &result);
            result->SetLength(batch_size);
            std::cout << "v * v = " << result;
            // Output: v * v = (1.0, 4.0, 9.0, ...
        \end{minted}
    \end{minipage}
    \hfill
    \begin{minipage}[t]{0.45\textwidth}
        \begin{minted}[fontsize=\small, frame=none]{julia}
            v = [1.0, 2.0, 3.0, 4.0]
            p = MakeCKKSPackedPlaintext(cc, v)
            c = Encrypt(cc, public_key, p)
            c_squared = EvalMult(cc, c, c)
            result = Plaintext()
            Decrypt(cc, private_key, c_squared, result)
            SetLength(result, batch_size)
            println("v * v = ", result)
            # Output: v * v = (1.0, 4.0, 9.0, ...
        \end{minted}
    \end{minipage}
    \caption{Ciphertext multiplication using OpenFHE in C++ (left) and OpenFHE.jl in Julia (right).}
    \label{lst:openfhe-vs-openfhejl}
\end{listing}

\begin{listing}[!htbp]
    \begin{minted}[frame=none]{julia}
        v = [1.0, 2.0, 3.0, 4.0]
        c = encrypt(v, public_key, context)
        c_squared = c * c
        result = decrypt(c_squared, private_key)
        println("v * v = ", result)
        # Output: v * v = [1.0, 4.0, 9.0, ...
    \end{minted}
    \caption{Ciphertext multiplication using SecureArithmetic.jl in Julia.}
    \label{lst:securearithmeticjl}
\end{listing}

\section{Accuracy and performance of the CKKS scheme}
\label{sec:accuracy_performance}

In this section, we analyze the accuracy and performance of basic CKKS operations as
provided by OpenFHE to create a baseline for the more complex numerical simulations
shown later.
We first describe the CKKS
configuration used throughout this paper and the method we use to measure errors and
runtime. We then present our accuracy and performance analysis for individual FHE
operations and discuss the implications of these results for the secure numerical
simulations. All code used for the experiments and our numerical results are
available in the reproducibility repository \cite{reproKholodSchlottkeLakemper24}.

\subsection{Experimental setup and measurement methodology}\label{sec:experimental_setup}

The experiments in this paper were conducted with the OpenFHE library v1.2.0 via the
SecureArithmetic.jl package, using a CKKS configuration that ensures 128-bit security
and uses $S_\text{word}$~=~64~bits as the native integer size. The security level is
determined based on the homomorphic encryption standard \cite{cryptoeprint:2019/939}, and
OpenFHE throws an error if any parameters provided violate the security guarantee. All
other parameters were chosen on a best-effort basis to balance performance and accuracy.
Unless noted otherwise, the configuration used for the OpenFHE setup is as shown in
Table~\ref{tab:configuration}.
While the OpenFHE library has many more user-configurable options, we only focus on those that need to be set explicitly or where we deviate from the default values.
For details on the CKKS configuration parameters, please
refer to the OpenFHE documentation \cite{OpenFHE}.
\begin{table}[htbp!]
    \centering
    \begin{tabular}{l l}
        \toprule
        \textbf{Parameter} & \textbf{Value} \\ \midrule
        Batch size & $2^5$--$2^{16}$ (\emph{as small as possible)} \\
        Bootstrapping level budget & $[4, 4]$ \\
        Enabled features & \texttt{ADVANCEDSHE},
        \texttt{FHE}, \texttt{KEYSWITCH}, \texttt{LEVELEDSHE}, \texttt{PKE}\\
        Size of first modulus $q_0$ & $60$ bits \\
        Size of scaling modulus $q_i$ & $59$ bits \\
        Scaling technique & \texttt{FLEXIBLEAUTO} \\
        Secret key distribution & \texttt{SPARSE\char`_TERNARY} (use
        \texttt{UNIFORM\char`_TERNARY} for production)\\
        Security level & \texttt{HEStd\char`_128\char`_classic} \\ \bottomrule
    \end{tabular}
    \caption{Parameters for the CKKS configuration in OpenFHE. For more details, please
    refer to the OpenFHE documentation \cite{OpenFHE}.}
    \label{tab:configuration}
\end{table}

For the batch size/capacity, we used the smallest power of two that still fits the data.
The multiplicative depth available after bootstrapping was chosen to
be $l_\text{refresh} = 15$, with a maximum multiplicative depth of $l_\text{max} = 33$.
The ring dimension is chosen automatically by OpenFHE, which usually corresponds to $N_R = 2^{17}$
given the remaining parameters. Since the security standards in the OpenFHE library are
being constantly improved, a different ring dimension may be supported
for these parameters in the future.
Furthermore, since the current homomorphic encryption standard \cite{cryptoeprint:2019/939}
does not provide CKKS parameters when using a sparse ternary distribution for
the secret key generation, in production settings it is advisable to use the uniform
ternary distribution.

The performance and accuracy investigations in this section require applying operations to
encrypted and unencrypted data. The main ciphertext is always initialized to
\begin{equation}
    \left[\sin\left(2 \pi \frac{1}{L}\right), \sin\left(2 \pi \frac{2}{L}\right), ..., \sin\left(2 \pi \frac{L}{L}\right)\right],
\end{equation}
with $L = 64$. Scalars are always initialized to $1 + \pi/30$, and plaintexts to
\begin{equation}
    \left[1+\frac{\pi}{30}, 1+\frac{\pi}{30}, ..., 1+\frac{\pi}{30}\right],
\end{equation}
also with a length of $L = 64$. Therefore, the batch size/capacity is always $64$.
If required, a second ciphertext was created by encrypting
the above plaintext.
Due to the approximate nature of CKKS, the encrypted values differ slightly from the
original.

Since the computational overhead of using OpenFHE through the Julia packages is negligible
compared to the cost of the FHE operations themselves, in this paper we use
SecureArithmetic.jl to analyze the performance.
The experiments
were conducted on an AMD Ryzen Threadripper 3990X 64 core processor with 256 GiB RAM using a
single thread.
All measurements discard the initial just-in-time compilation time of the Julia programming
language. We measured the runtime of FHE operations by averaging over five consecutive
executions.

Numerical errors are measured by comparing the decrypted result of a secure operation with a
corresponding unencrypted operation. In case of encryption, decryption, and bootstrapping
operations, the error is calculated relative to doing nothing with the unencrypted data.
Secure addition, multiplication, and rotation are compared to regular addition,
multiplication, and Julia's \texttt{Base.circshift}. For the error measurements, we used the
$L^\infty$ norm, 
which denotes the maximum absolute element-wise error.
The observed error depends heavily on the cryptographic parameters such as
ring dimension, security level, or scaling modulus. Therefore, error measurements should not
be interpreted by their absolute numbers, but rather by their order of magnitude and trend.

\subsection{Plaintext encoding and decoding}
\begin{figure}[!htbp]
    \centering
    \begin{subfigure}[b]{0.48\textwidth}
        \centering
        \includegraphics[width=\textwidth]{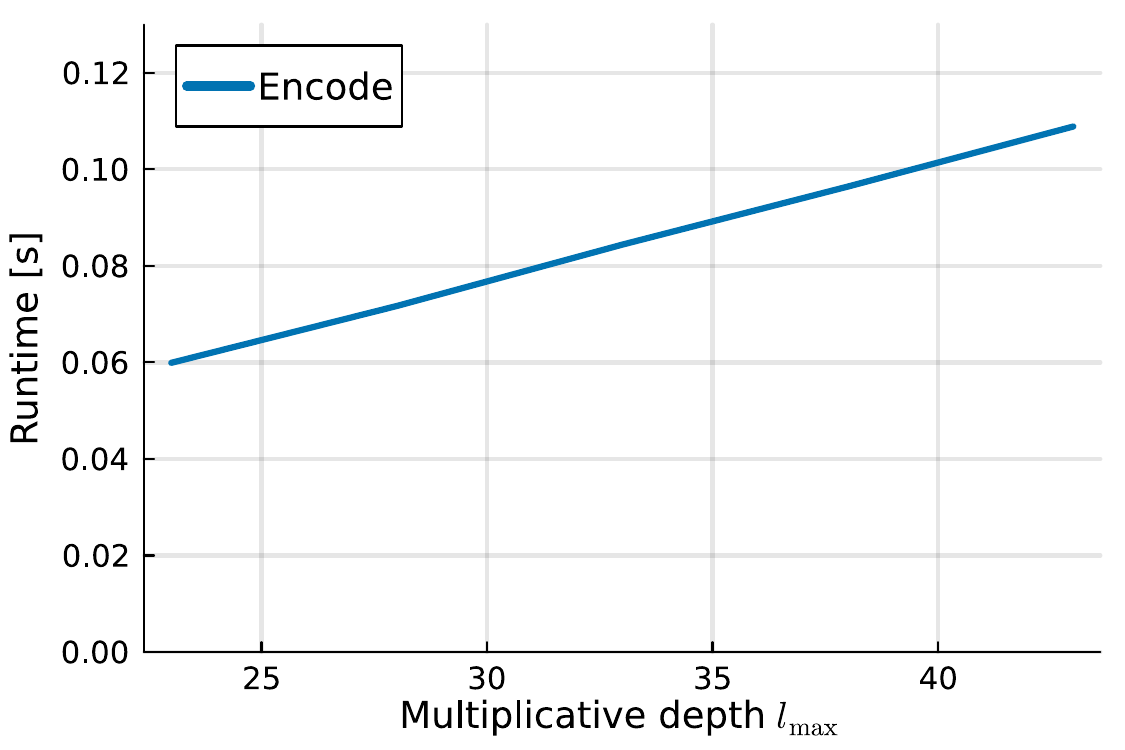}
        \label{fig:performance-encoding-decoding}
    \end{subfigure}
    \caption{Performance evaluation for encoding a vector of real numbers into a plaintext.}
    \label{fig:encoding-decoding}
\end{figure}

As described in Sec.~\ref{sec:encoding}, the encoding operation is used to
transform a vector of real numbers into a plaintext, where it is stored as a polynomial. The
decoding is the inverse operation, which converts the polynomial back into a data vector.
Unfortunately, the OpenFHE library does not provide a direct way to decode a plaintext (as extra noise is added during decoding for security purposes), thus
we can only measure the runtime of the encoding operation, but not the decoding operation, nor
the error of either operation. From Fig.~\ref{fig:encoding-decoding}, we can see that
runtime required for encoding a data vector into a plaintext grows linearly with the
multiplicative depth $l_\text{max}$.
While generally the plaintext encoding is independent of $l_\text{max}$, by default
the OpenFHE library already prepares the plaintext for multiplication with a
ciphertext, i.e., transforms the plaintext to ciphertext space. This
requires additional expensive operations, which depend on
$l_\text{max}$, to be performed during encoding.
The level used to prepare the plaintext can be configured (changed from its default value) by the user of the library.

\subsection{Encryption and decryption}
We proceed with analyzing the accuracy and performance of the two most fundamental
operations of cryptography: encryption and decryption. Since it is impossible to measure the
error of these operations independently, the error was determined after encrypting and
immediately decrypting a data vector.

\begin{figure}[!htbp]
    \centering
    \begin{subfigure}[b]{0.48\textwidth}
        \centering
        \includegraphics[width=\textwidth]{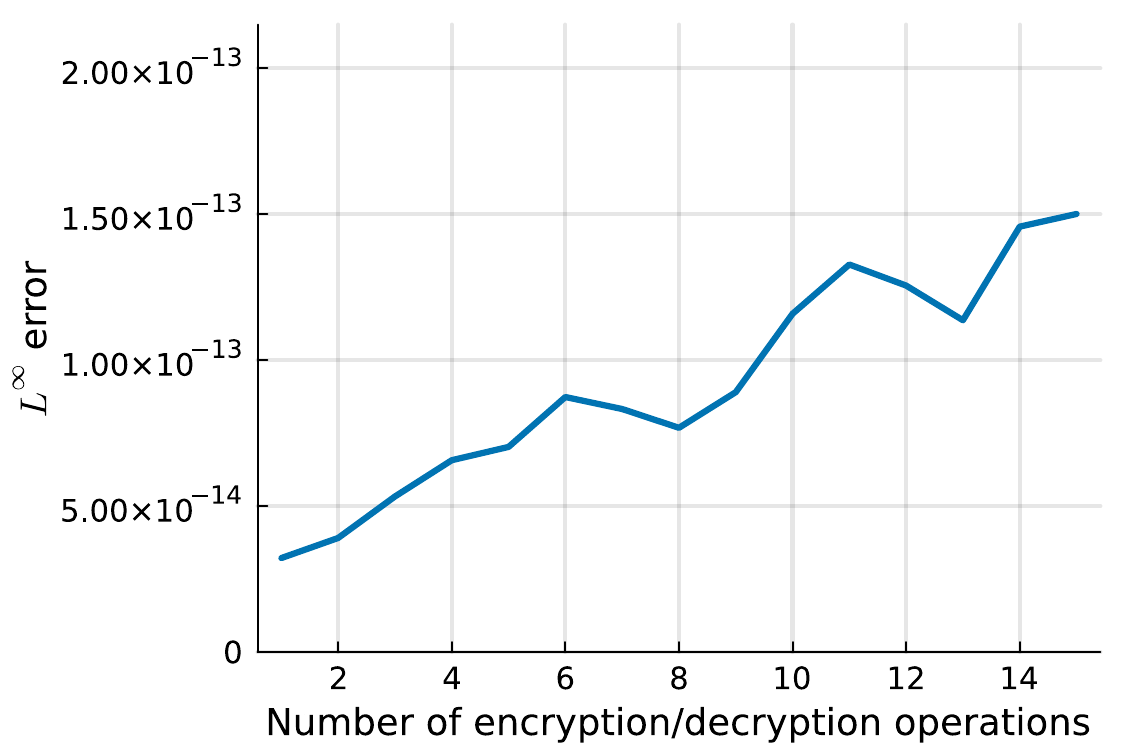}
        \caption{Error from encryption/decryption operations.}
        \label{fig:accuracy-encryption-decryption}
    \end{subfigure}
    \hfill
    \begin{subfigure}[b]{0.48\textwidth}
        \centering
        \includegraphics[width=\textwidth]{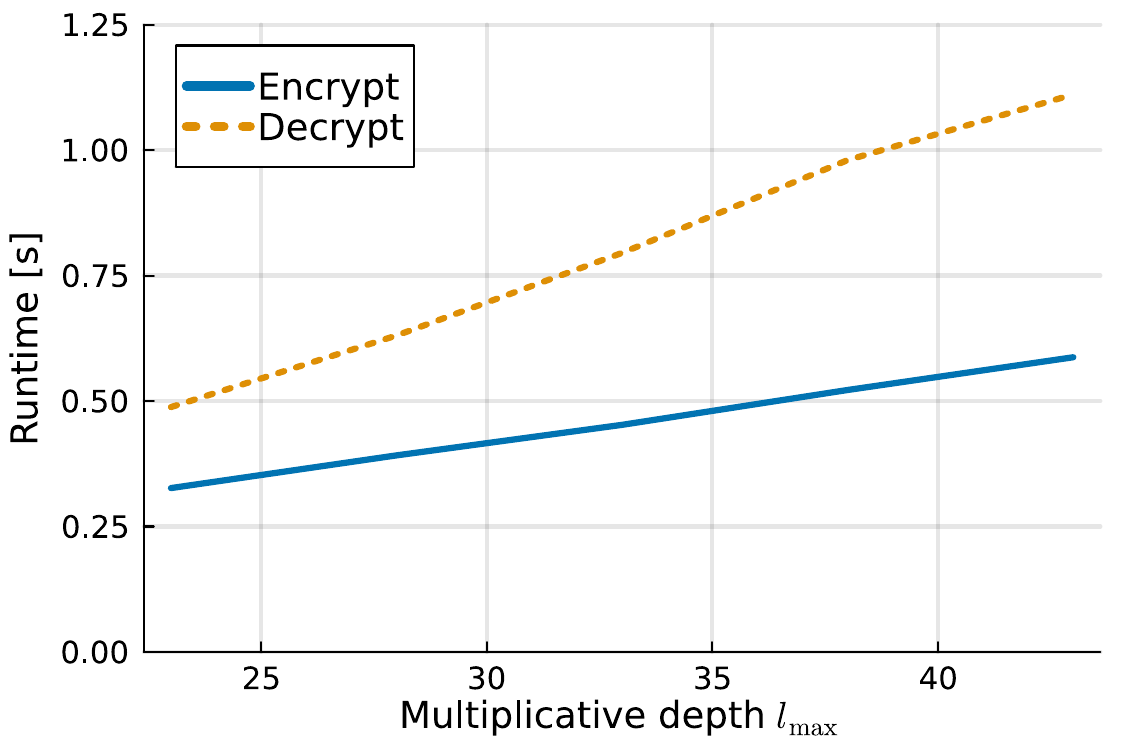}
        \caption{Time required for encryption/decryption operations.}
        \label{fig:performance-encryption-decryption}
    \end{subfigure}
    \caption{Error analysis (left) and performance evaluation (right) for encryption and
    decryption of a ciphertext. For the error investigation, the multiplicative depth was
    set to $l_\text{max} = 33$.}
    \label{fig:encryption-decryption}
\end{figure}

As shown in Fig.~\ref{fig:accuracy-encryption-decryption}, the error of a single
encryption/decryption operation is of $\mathcal{O}(10^{-13})$, and remains at that level
even after more than ten subsequent encryption/decryption pairs.
Since the error is
approximately in the range of machine precision (for the parameters summarized in Table~\ref{tab:configuration}), and usually only a single encryption and
decryption are performed, the error is negligible for all practical purposes.

As mentioned in Section~\ref{sec:levels}, we expect the runtime of ciphertext operations to depend linearly on the
multiplicative depth $l_\text{max}$. This is confirmed by the results in
Fig.~\ref{fig:performance-encryption-decryption}. It is also clear that encryption and
decryption take a non-negligible amount of time, already in the order of a second.

Due to the nature of how we measure the accuracy of other operations, the
error incurred by one encryption/decryption operation will be included in the error analysis
of all other FHE operations in the subsequent sections. Moreover, since it is impossible
to measure the error of encoding and decoding operations independently, they are also implicitly included.
A more detailed discussion of encryption/decryption errors can be found in
\cite{KimPapadimitriouEtAl20}.

\subsection{Addition}
\begin{figure}[!htbp]
    \centering
    \begin{subfigure}[b]{0.48\textwidth}
        \centering
        \includegraphics[width=\textwidth]{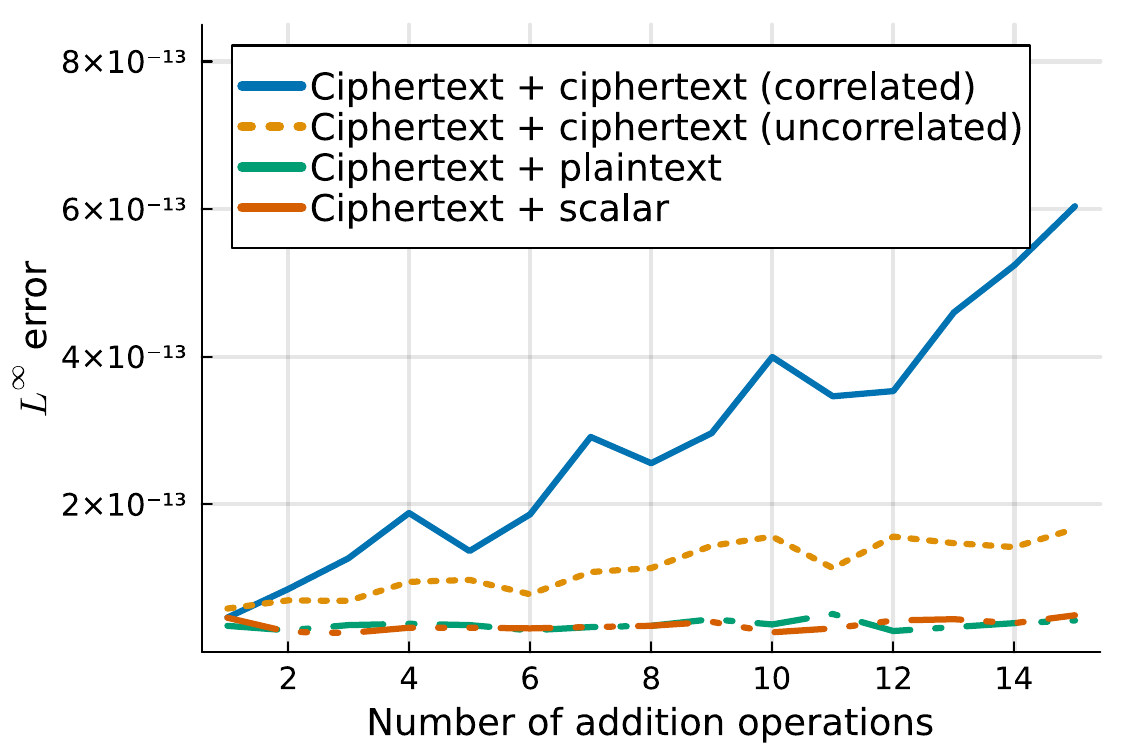}
        \caption{Error after addition operations.}
        \label{fig:accuracy-addition}
    \end{subfigure}
    \hfill
    \begin{subfigure}[b]{0.48\textwidth}
        \centering
        \includegraphics[width=\textwidth]{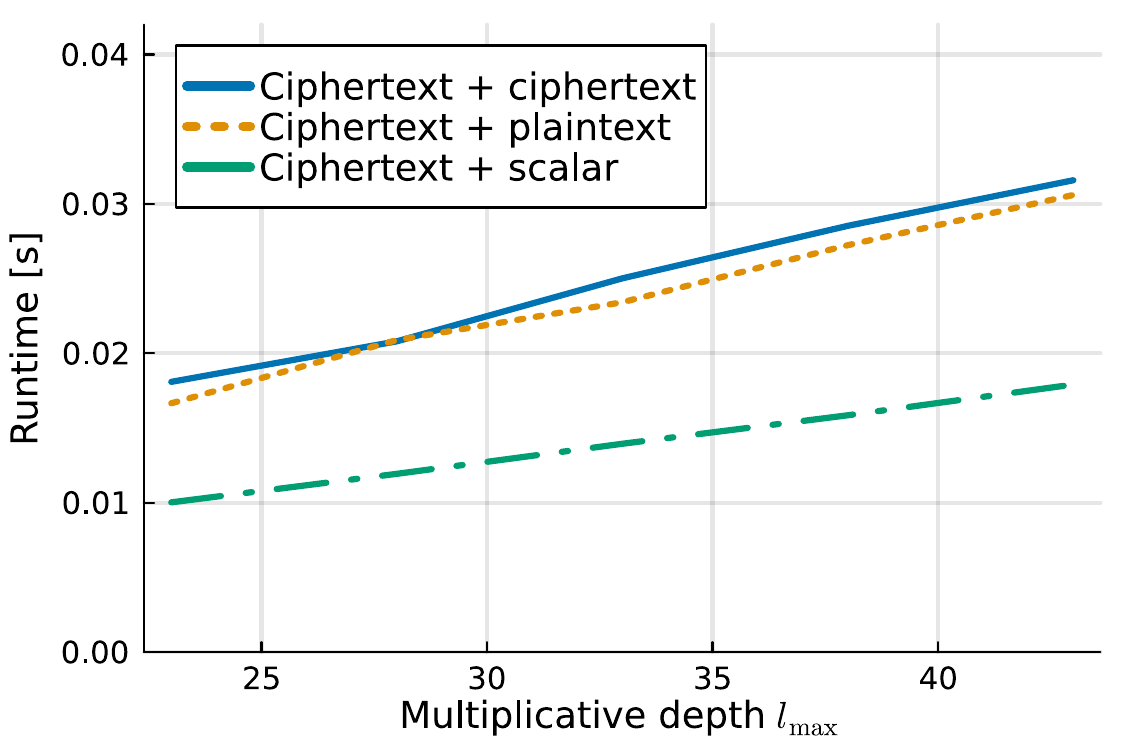}
        \caption{Time required for addition operations.}
        \label{fig:performance-addition}
    \end{subfigure}
    \caption{Error analysis (left) and performance evaluation (right) for addition
    operations between ciphertexts, ciphertexts and plaintexts, and ciphertexts and scalars.
    For the error investigation, the multiplicative depth was set to $l_\text{max} = 33$.}
    \label{fig:addition}
\end{figure}

Error and runtime measurements for addition operations are presented in
Fig.~\ref{fig:addition}. Here, we analyze the addition of two ciphertexts, a ciphertext and
a plaintext, and a ciphertext and a scalar. The error of adding a plaintext or a scalar
value to a ciphertext is virtually negligible at $\mathcal{O}(10^{-14})$, which is
approximately the truncation error for double-precision arithmetic. Adding two ciphertexts
incurs an error only slightly larger, on the order of $\mathcal{O}(10^{-13})$
(Fig.~\ref{fig:accuracy-addition}).

When adding ciphertexts, it is instructive to distinguish between ciphertexts with
correlated and uncorrelated errors.
Two ciphertexts have correlated errors if they share the same noise or if their noise terms
are not statistically independent. This typically happens when the same ciphertext is used
more than once in an algorithm, e.g., when adding a ciphertext to itself multiple times.
To achieve uncorrelated ciphertext errors, one needs to encrypt the second ciphertext
summand independently for each operation (or rerandomize the ciphertext by adding an
encryption of zero). As shown in
Fig.~\ref{fig:accuracy-addition}, the difference between these two cases is significant.
Although the errors initially have the same magnitude, correlated errors increase linearly
with the number of additions, as the same errors are repeatedly added. In contrast,
uncorrelated errors grow at a rate proportional to the square root of the number of
additions, matching the expected behavior based on the Central Limit Theorem.

As with the other operations before, the runtime of the addition operation grows linearly
with the multiplicative depth $l_\text{max}$. While the addition of two ciphertexts or a
ciphertext and a plaintext requires a similar amount of time, adding a scalar to a ciphertext
is considerably faster.

\subsection{Multiplication}
\begin{figure}[!htbp]
    \centering
    \begin{subfigure}[b]{0.48\textwidth}
        \centering
        \includegraphics[width=\textwidth]{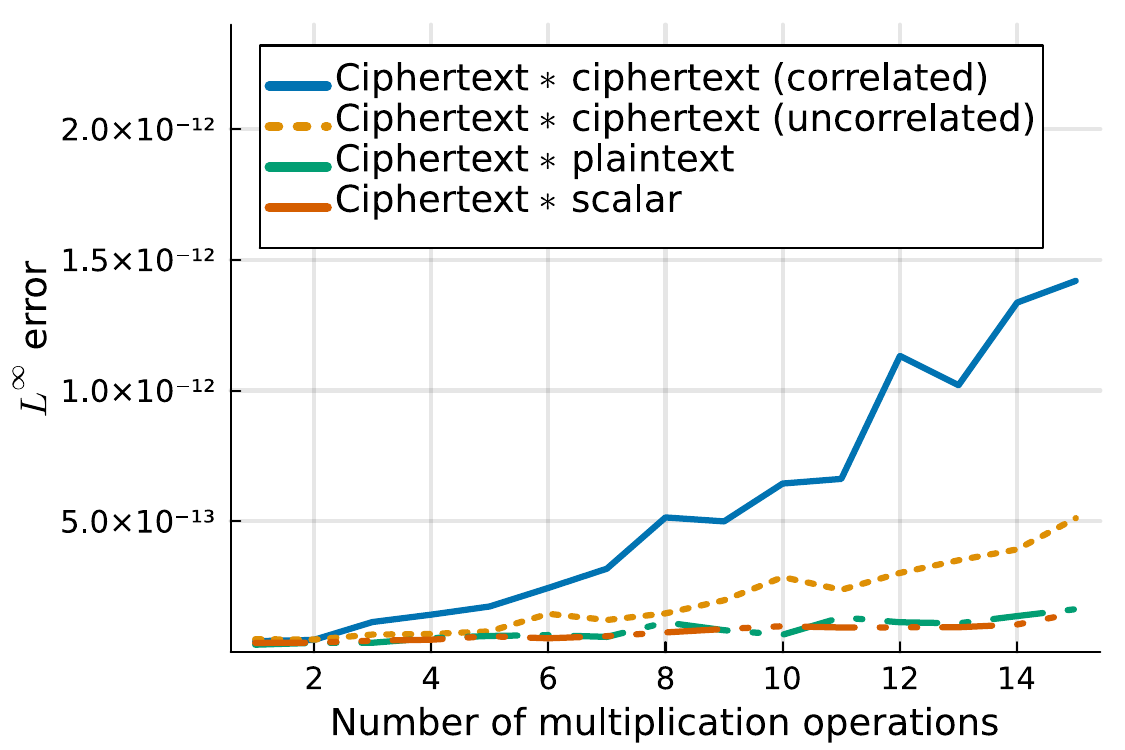}
        \caption{Error after multiplication operations.}
        \label{fig:accuracy-multiplication}
    \end{subfigure}
    \hfill
    \begin{subfigure}[b]{0.48\textwidth}
        \centering
        \includegraphics[width=\textwidth]{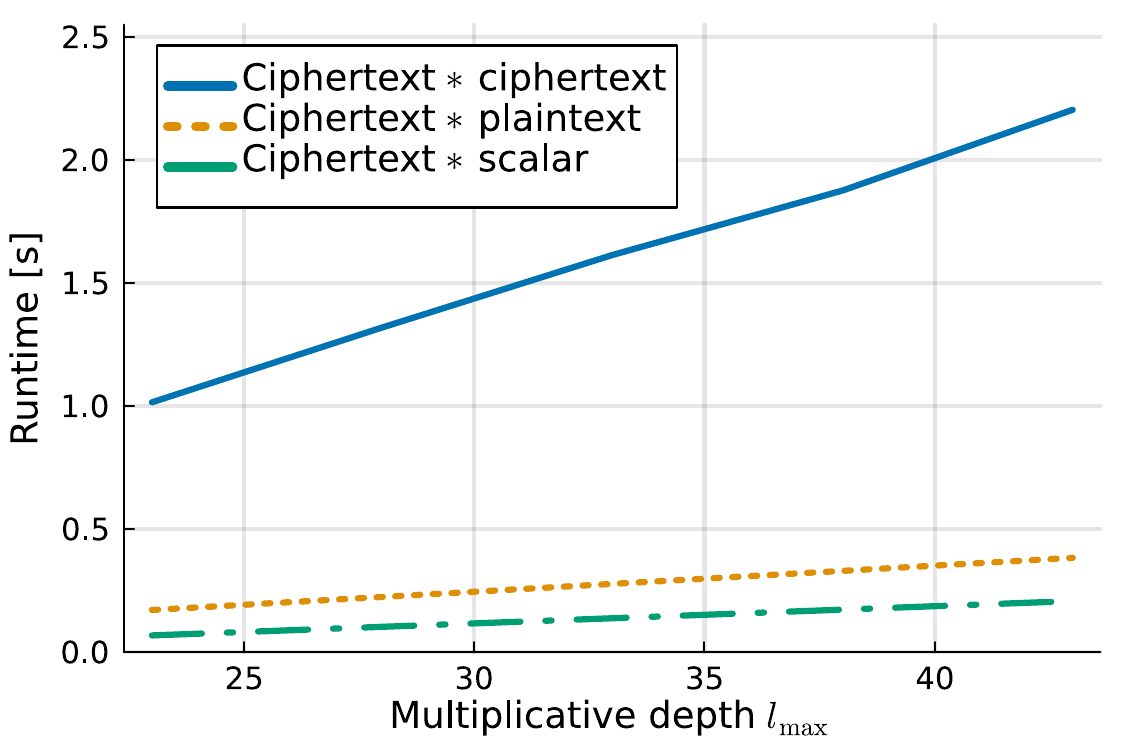}
        \caption{Time required for multiplication operations.}
        \label{fig:performance-multiplication}
    \end{subfigure}
    \caption{Error analysis (left) and performance evaluation (right) for multiplication
    operations between ciphertexts, ciphertexts and plaintexts, and ciphertexts and scalars.
    For the error investigation, the multiplicative depth was set to $l_\text{max} = 33$.}
    \label{fig:multiplication}
\end{figure}

Next, we analyze the accuracy and performance of multiplication operations.
As shown in Fig.~\ref{fig:accuracy-multiplication}, the behavior of the error from
multiplying two ciphertexts depends on whether the ciphertexts have correlated or
uncorrelated errors.
Similar to the findings for ciphertext addition, multiplying ciphertexts with correlated
errors leads to a linear increase in the overall error, while multiplying ciphertexts with
uncorrelated errors results in an error growth proportional to the square root of the number
of operations. In both cases, we adjusted the levels of the ciphertexts used for
multiplication to maintain consistency across the ciphertexts. In practice, the total error
from multiplication in an FHE algorithm will be between the correlated and uncorrelated error
curves.

Unlike addition or multiplication of a ciphertext with a plaintext or a scalar,
multiplication of two ciphertexts requires a relinearization step, which is especially
computationally expensive due to the necessary key switching operation. Therefore, the
runtime of ciphertext-ciphertext multiplication increases significantly compared to
addition, as can be observed in Fig.~\ref{fig:performance-multiplication}, and is much
larger than for multiplication by plaintext or scalar. Furthermore, multiplication of a
ciphertext by a plaintext or scalar is much slower than the corresponding addition
operations due to the necessary rescaling operation. Generally, one should keep in mind that
the runtime for key switching is always significantly larger than the runtime of rescaling,
which in turn is significantly larger than the runtimes of the underlying component-wise
multiplication or addition.

As before, the runtime grows with
increased multiplicative depth $l_\text{max}$. For the performance measurements,
we evaluated the runtime of the second multiplication of the ciphertext with another factor,
since in OpenFHE with the \texttt{FLEXIBLEAUTO} rescaling technique, the rescaling step is
not applied to the first multiplication.

\subsection{Rotation}
\begin{figure}[!htbp]
    \centering
    \begin{subfigure}[b]{0.48\textwidth}
        \centering
        \includegraphics[width=\textwidth]{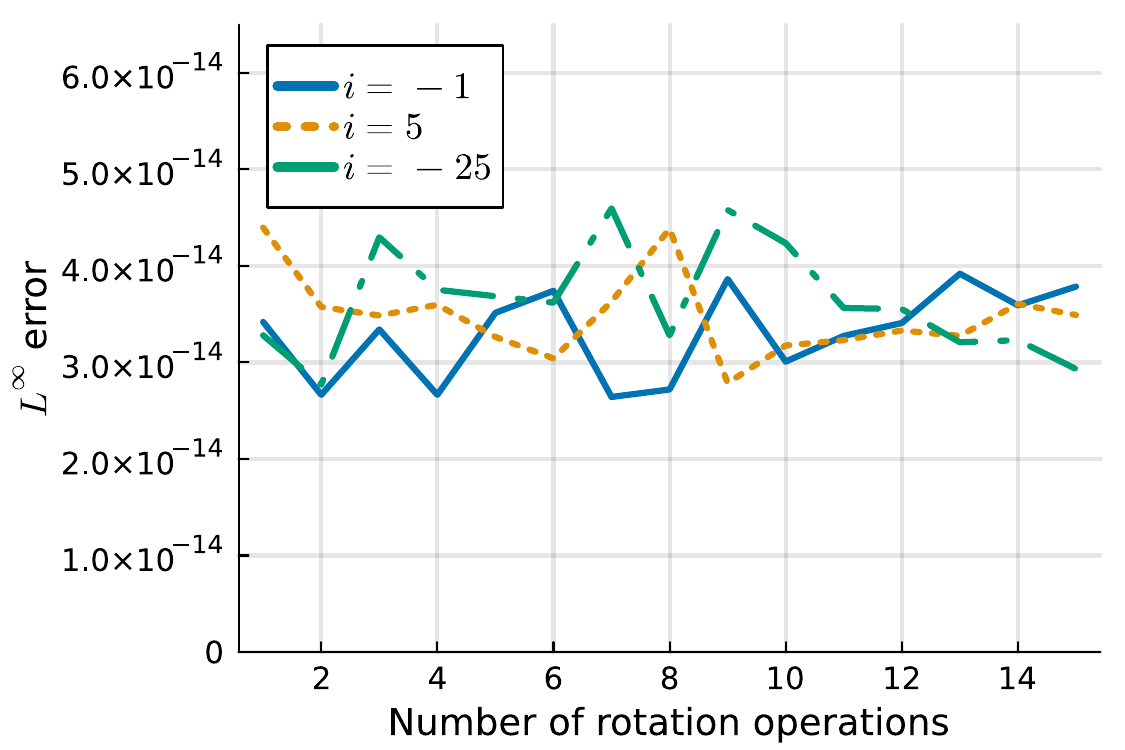}
        \caption{Error after rotation operations.}
        \label{fig:accuracy-rotation}
    \end{subfigure}
    \hfill
    \begin{subfigure}[b]{0.48\textwidth}
        \centering
        \includegraphics[width=\textwidth]{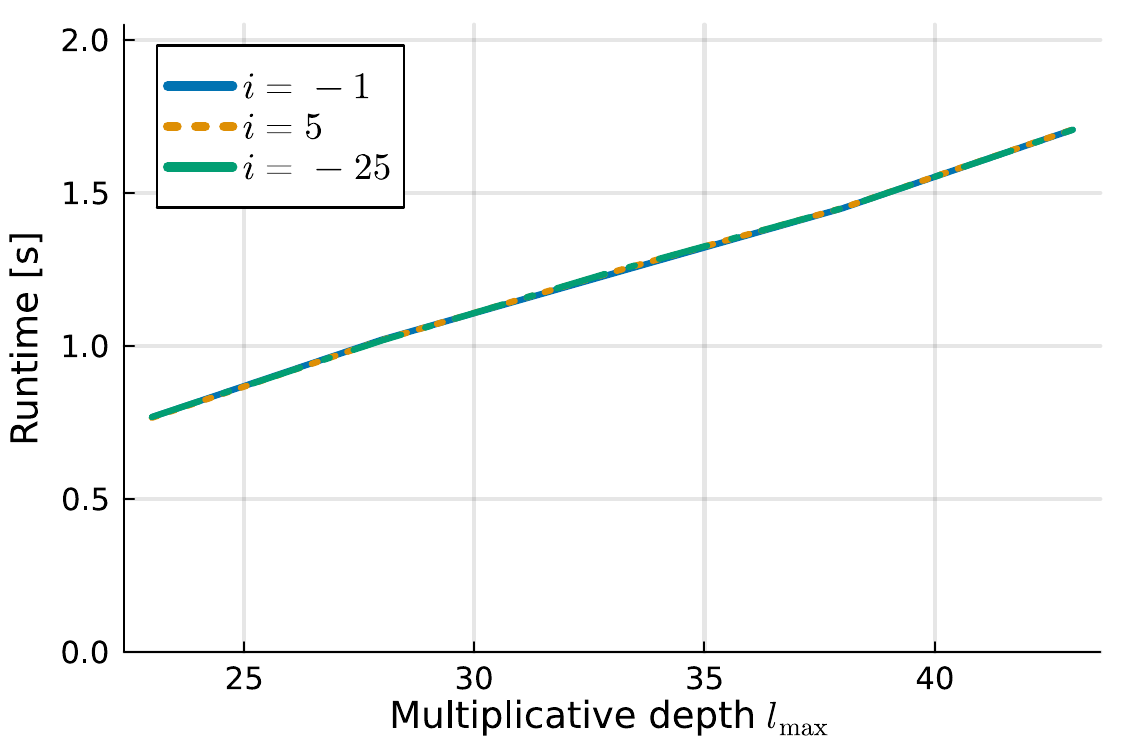}
        \caption{Time required for rotation operations.}
        \label{fig:performance-rotation}
    \end{subfigure}
    \caption{Error analysis (left) and performance evaluation (right) for rotation
    operations by different rotation indices $i \in \{-1, 5, -25\}$. For the error investigation, the multiplicative depth was set to $l_\text{max} = 33$.}
    \label{fig:rotation}
\end{figure}

The impact of rotation operations on error and execution time is shown in
Fig.~\ref{fig:rotation} for different rotation indices. The rotation operation does not
result in error growth, thus the error remains at the level of a single
encryption-decryption operation. In terms of the runtime, there is still a linear growth
with the multiplicative depth $l_\text{max}$. Neither the error nor the runtime depend on
the rotation index.

\subsection{Bootstrapping}\label{sec:accuracy-performance-bootstrapping}
Finally, we analyze the accuracy and performance of the bootstrapping procedure. The errors
are measured without any arithmetic operations between two subsequent bootstrapping
operations. To make the results comparable between different settings for the maximum
multiplicative depth $l_\text{max}$, we always measure the execution time after dropping
the ciphertext to two remaining multiplicative levels $l = 2$ before performing the
bootstrapping operation.

To improve the accuracy of the bootstrapping operation, in \cite{bae2022metabts}
an iterative (or multiprecision) bootstrapping technique was introduced. The idea is to
perform multiple bootstrapping operations in sequence when refreshing the ciphertext, each time progressively reducing the error by the precision (in bits) of CKKS bootstrapping. This
can significantly increase the bits of precision available for an FHE computation, as
demonstrated in, e.g., \cite{badawi2023demystifying}. The OpenFHE library supports iterative
bootstrapping with two consecutive iterations, which requires an experimentally
determined precision of a single bootstrapping as an input. In this paper, we use both the
standard bootstrapping method with a single iteration and the iterative procedure with two
iterations, setting the experimentally determined precision to $19$ bits.

\begin{figure}[!htbp]
    \centering
    \begin{subfigure}[b]{0.48\textwidth}
        \centering
        \includegraphics[width=\textwidth]{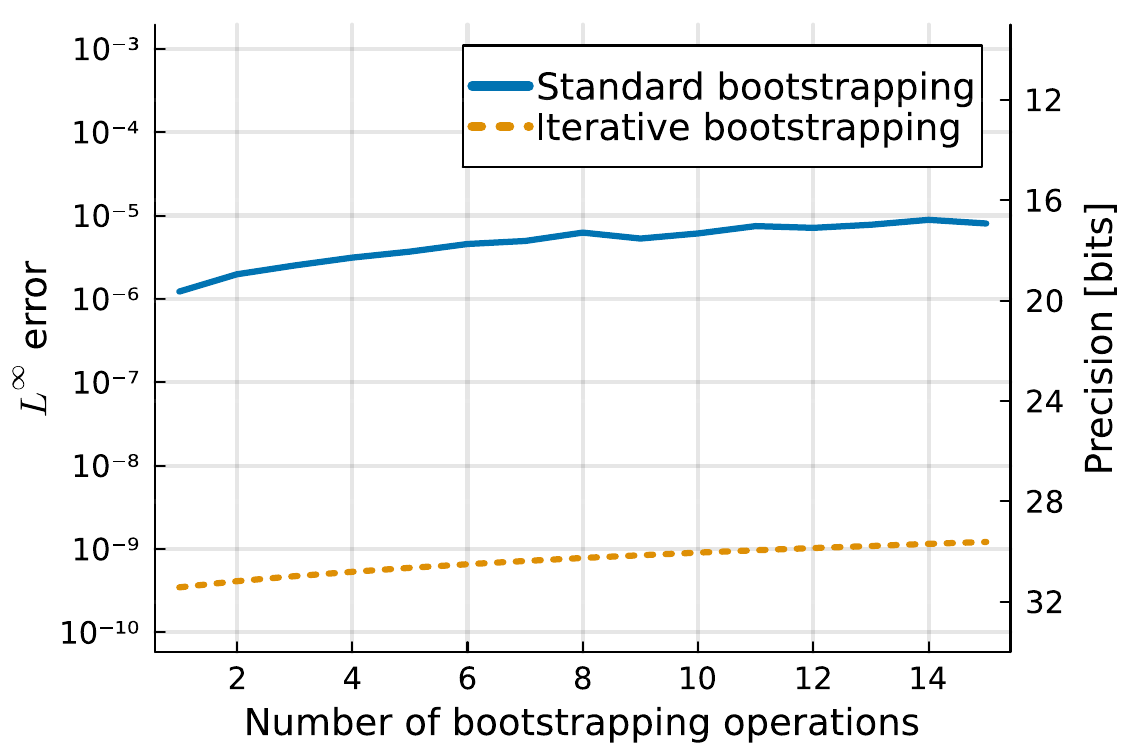}
        \caption{Error after bootstrapping operations.}
        \label{fig:accuracy-bootstrapping}
    \end{subfigure}
    \hfill
    \begin{subfigure}[b]{0.48\textwidth}
        \centering
        \includegraphics[width=\textwidth]{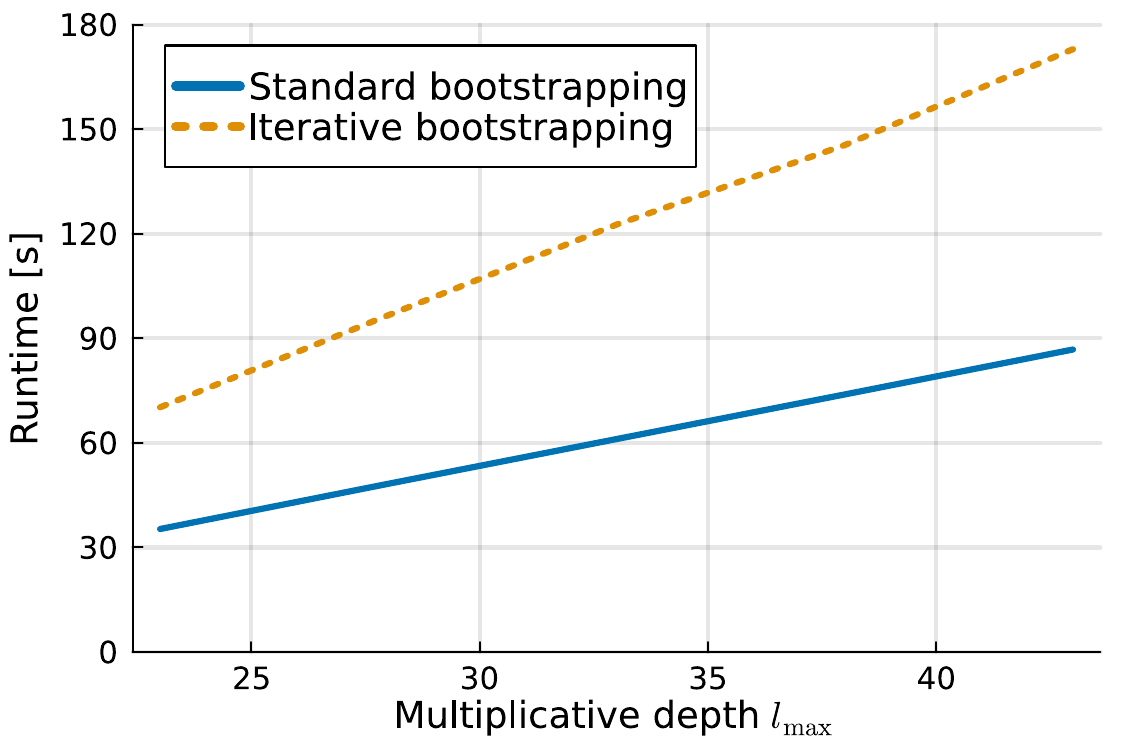}
        \caption{Time required for bootstrapping operations.}
        \label{fig:performance-bootstrapping}
    \end{subfigure}
    \caption{Error analysis (left) and performance evaluation (right) for bootstrapping
    operations. For the error investigation, the multiplicative depth was set to
    $l_\text{max} = 33$. The standard bootstrapping procedure uses a single bootstrapping
    operation, while the iterative procedure performs two subsequent bootstrapping
    operations.}
    \label{fig:bootstrapping}
\end{figure}

From Fig.~\ref{fig:accuracy-bootstrapping}, it becomes clear that bootstrapping has a much higher approximation error than all other CKKS operations, which is consistent with previous findings,
e.g., in \cite{Bootstrapping}.
The $L^\infty$ error at $\mathcal{O}(10^{-6})$ for standard bootstrapping and
$\mathcal{O}(10^{-9})$ for iterative bootstrapping is orders of magnitude greater than in
all previous experiments (adding further bootstrapping iterations is expected to
make this error comparable to the error after leveled CKKS operations).
When performing multiple bootstrapping procedures in a row, a gradual decline in precision
is observed, with the rate of precision loss becoming progressively smaller. This is
consistent with~\cite{Bootstrapping}, where the error eventually reaches quasi-steady state
conditions after many bootstrapping invocations. The primary cause for this behavior is the
approximation of the modular reduction---the function evaluated as part of CKKS
bootstrapping---with a sine wave. This approximation is only accurate near zero within a
small, periodically repeating interval, and its error grows as the message value moves
farther from zero~\cite{Bootstrapping}. Each application of the sine wave approximation
introduces additional error, but the precision loss decreases with each step, as each
bootstrapping operates on already-degraded input (rather than the more accurate value fed
into the first bootstrapping). In practice, this means that CKKS can be used for very
deep computations with many hundreds of bootstrapping invocations, with a small cumulative
loss in precision compared to the first bootstrapping. Moreover, CKKS also provides tools to
further reduce the noise using Hermite
interpolations~\cite{alexandru2024functionalbootstrapping}, thereby enabling practically
unlimited deep computations.

Fig.~\ref{fig:performance-bootstrapping}
further reveals that bootstrapping is also the most time-consuming of all CKKS operations,
with a runtime that again grows linearly with the maximum multiplicative depth
$l_\text{max}$. The runtime for the iterative bootstrapping is approximately twice as large
as for the standard bootstrapping. Thus in practice, the increased accuracy of the iterative
bootstrapping technique has to be balanced with the additional computational cost.
Furthermore, we have to emphasize that the errors and runtimes reported here are highly
dependent on the specific CKKS setup and the actual user data, and should thus be only
interpreted as indicators and not taken for their absolute values. Unless noted otherwise,
in the remainder of the paper we use the standard bootstrapping method.

\begin{figure}[!htbp]
    \centering
    \includegraphics[width=0.8\textwidth, trim=8mm 2mm 10mm 2mm, clip]{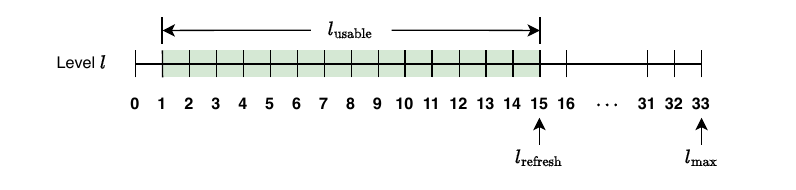}
    \caption{Overview of an exemplary CKKS bootstrapping process with maximum
    multiplicative depth $l_\text{max} = 33$. Standard bootstrapping requires one level to
    get started, and then consumes $18$ levels during the actual process. After
    bootstrapping, the ciphertext level is restored to $l_\text{refresh} = 15$ levels, and
    $l_\text{usable} = 14$ levels may be used before the next bootstrapping.}
    \label{fig:bootstrapping-overview}
\end{figure}

At this point it should be mentioned that the bootstrapping procedure (as implemented
in OpenFHE) requires one remaining multiplicative level per bootstrapping iteration to get
started and also consumes multiple levels for the bootstrapping itself.
Therefore, after bootstrapping the ciphertext level is restored to $l = l_\text{refresh} <
l_\text{max}$, and the number of levels for continuous operation is $l_\text{usable} =
l_\text{refresh} -1$ for standard and $l_\text{usable} = l_\text{refresh} -2$ for iterative
bootstrapping. Fig.~\ref{fig:bootstrapping-overview} illustrates these values, which are
important to consider when designing an FHE application.
A more detailed discussion of the slightly different behavior of the corresponding OpenFHE
function \texttt{GetLevel()} is provided in \ref{app:openfhe_details}.

We further observed in our experiments that in addition to $l_\text{max}$, the error and
runtime of bootstrapping are also sensitive to many other parameters, such as the ring
dimension $N_R$, the batch size/capacity, the ciphertext modulus and scaling factor, the
bootstrapping level budget, and the multiplicative depth available after bootstrapping
$l_\text{refresh}$. It is thus difficult to quantify the accuracy and performance behavior
of bootstrapping in general. Moreover, the error also depends on the absolute value of the
data in the ciphertext. For algorithms that use bootstrapping, it is thus recommended to
rescale the data such that its magnitude is less than one \cite{CheonHanEtAl18}. While it is
beyond the scope of this paper to analyze the performance and error implications of all
these parameters, we recommend to carefully tailor them to the specific application.

To summarize, of all CKKS operations, bootstrapping is the most detrimental to the accuracy
and performance of any FHE application, and should be used as sparingly as possible.
However, it is also the ingredient that enables an unlimited number of multiplications in an
FHE algorithm (for most of practical applications) and thus facilitates the evaluation of
truly arbitrary functions on encrypted data. When the accuracy of the computations is a
primary concern, iterative bootstrapping with two iterations can reduce the error by multiple
orders of magnitude at the cost of doubling the runtime.

\subsection{Overview of accuracy and performance of OpenFHE operations}
As noted at the beginning of the accuracy and performance investigation, the exact values of
the errors and execution times depend on many factors such as the chosen CKKS parameters,
the actual operations performed, the data itself, and on the used compute system.
Consequently, all absolute values presented here should be regarded as indicative estimates.
We summarize our findings in Table~\ref{tab:operations}, where we analyze the errors and
execution times of different FHE operations relative to ciphertext-ciphertext addition. From
the table, it is again evident that the most inaccurate and time-consuming FHE operation is
bootstrapping. Hence, it is important to use bootstrapping as rarely as possible.
At the same time, the previous results also show a linear dependence of the runtime on
the maximum multiplicative depth $l_\text{max}$ for all FHE operations.
Although increasing
$l_\text{max}$ and the usable levels after bootstrapping $l_\text{usable}$ reduces the
number of required bootstrapping operations, the runtime of all operations (including
bootstrapping itself) increases linearly. Therefore, $l_\text{max}$ should be chosen
carefully to balance error and runtime, and its optimal value depends on the specific
application.

\begin{table}[!htbp]
    \centering
    \begin{tabular}{l l l}
        \toprule
        Operation & Relative error & Relative runtime\\
        \midrule
        Encrypt plaintext & $<1$ & $\approx 20$\\[1ex]
        Decrypt ciphertext & $<1$ & $\approx 30$\\
        \midrule
        Ciphertext + ciphertext & $ = 1$ & $ = 1$\\[1ex]
        Ciphertext + plaintext & $< 1$ & $\approx 1$\\[1ex]
        Ciphertext + scalar & $< 1$ & $< 1$\\
        \midrule
        Ciphertext $*$ ciphertext & $\approx 1$ & $\approx 60$\\[1ex]
        Ciphertext $*$ plaintext & $< 1$ & $\approx 10$\\[1ex]
        Ciphertext $*$ scalar & $< 1$ & $\approx 6$\\
        \midrule
        Rotate ciphertext & $< 1$ & $\approx 50$\\
        \midrule
        Bootstrap ciphertext (standard)   & $\approx 30,\!000,\!000$ & $\approx 2,\!500$\\[1ex]
        Bootstrap ciphertext (iterative) & $\approx 5,\!000$ & $\approx 5,\!000$\\
        \bottomrule
    \end{tabular}
    \caption{Error and runtime of OpenFHE operations relative to the
    addition of two ciphertexts, rounded to the first significant digit. The absolute
    $L^\infty$ error and runtime for ciphertext-ciphertext addition are
    $\mathcal{O}(10^{-14})$ and $\mathcal{O}(10^{-2})$~s, respectively. Independently
    encrypted ciphertexts were used for the binary ciphertext operations.}
    \label{tab:operations}
\end{table}

There are some additional measures that can be taken to improve the accuracy and performance
of the CKKS scheme. To reduce the error, it is recommended to increase the scaling factor. The iterative bootstrapping method decreases the errors
by several orders of magnitude compared to the standard bootstrapping method. Furthermore,
the OpenFHE library allows one to set $S_\text{word} = 128$~bits to internally use larger
integers. Most of these measures will, however, increase the runtime of the CKKS operations.
To improve the performance, it is advisable to always choose the smallest possible batch
size/capacity, since larger batch sizes increase the bootstrapping runtime and reduce the precision. Finally, it is
important to find the optimal multiplicative depth $l_\text{max}$: a larger one may reduce
the frequency of bootstrapping but at the same time increase the size of the ciphertext and,
thereby, the execution time of all FHE operations.

\section{Setting up secure numerical simulations}\label{sec:secure_simulation}

The central goal of this work is to demonstrate that fully homomorphic encryption can be used
for secure numerical simulations of partial differential equations (PDEs).
Since the CKKS scheme only supports three basic arithmetic operations (addition,
multiplication, and rotation), we need to recast the
numerical methods into a form that can be represented with the available operations. A
secondary goal for this manuscript is therefore to start developing the necessary
algorithmic building blocks, which will also be useful for other scientists who want to use
FHE for their own numerical applications.

In this section, we begin by introducing the linear advection equation in one and two
spatial dimensions, which will serve as the prototypical PDEs for our secure numerical
simulations. We then present two finite difference schemes which we use to
discretize these equations. Next, we rewrite these schemes in terms of the FHE
primitives provided by the CKKS scheme. Finally, we provide the full algorithm for
secure numerical simulations using FHE.

\subsection{Linear scalar advection equations}\label{sec:advection}

The linear scalar advection equation describes the transport of a
scalar field $u$ at constant speed. In one dimension, it is given by
\begin{equation}
    \pdv{u}{t} + a_x\pdv{u}{x} = 0,
    \label{AD1}
\end{equation}
with $u = u(t, x)$, where $t$ is time, $x$ is the spatial coordinate, and
$a_x > 0$ is the advection speed. In two dimensions, the equation extends to
\begin{equation}
    \pdv{u}{t} + a_x\pdv{u}{x} + a_y\pdv{u}{y} = 0,
    \label{AD2}
\end{equation}
with $u = u(t, x, y)$, and positive speeds $a_x, a_y > 0$.

Simulations are conducted on the domains $\Omega = [0, 1]$ in 1D and
$\Omega = [0, 1] \times [0, 1]$ in 2D, using periodic boundary conditions. Initial conditions are specified in the respective examples.

\subsection{Finite difference discretization}\label{sec:finite_difference}

\begin{figure}[!htbp]
    \centering
    \begin{subfigure}[t]{0.48\textwidth}
        \centering
	    \includegraphics[width=0.6\textwidth]{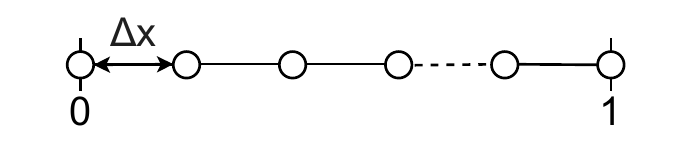}
        \caption{1D mesh with equidistant nodes on $\Omega = [0, 1]$.}
        \label{fig:mesh1d}
    \end{subfigure}
    \hfill
    \begin{subfigure}[t]{0.48\textwidth}
        \centering
        \includegraphics[width=0.6\textwidth]{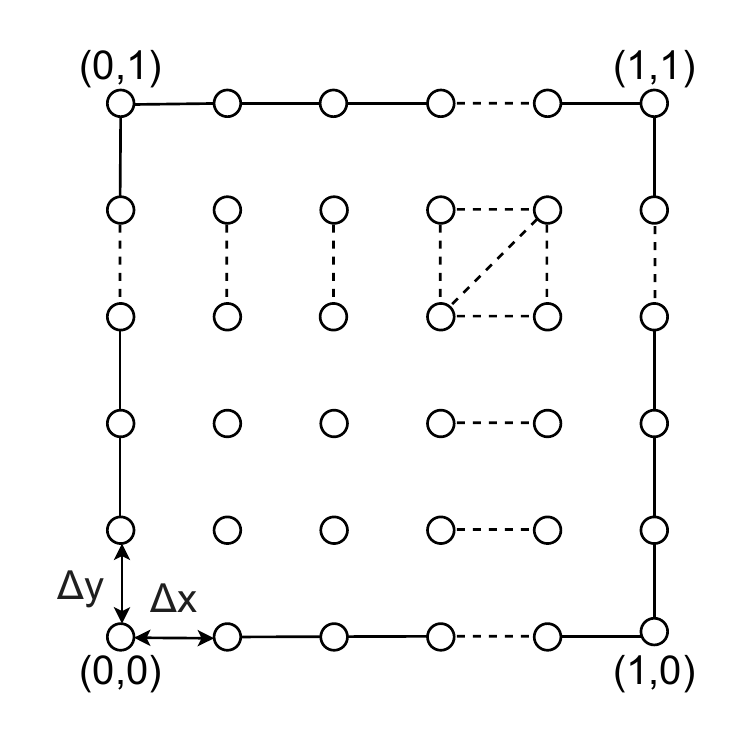}
        \caption{2D mesh with equidistant nodes on $\Omega = [0, 1] \times [0, 1]$.}.
        \label{fig:mesh2d}
    \end{subfigure}
    \caption{Cartesian meshes with equidistant node distributions for the finite-difference discretizations.}
    \label{fig:mesh}
\end{figure}

We use finite difference schemes to approximately solve the linear advection equations in
space and time. The computational domain is discretized by a Cartesian mesh with
equidistant nodes in each spatial direction as shown in Fig.~\ref{fig:mesh1d} for 1D and in
Fig.~\ref{fig:mesh2d} for 2D.

The numerical solution $u(t, x)$ is represented in 1D by the values $u_i^n = u(t_n, x_i)$ at
the mesh node locations $x_i = (i - 1) \Delta x$, with $i = 1, \dots, N_x$, and at time $t_n
= n \Delta t$, with $n = 0, 1, 2, \dots$ up to the final time. $\Delta x$ is the distance between two neighboring
nodes in the $x$-direction and $\Delta t$ is the time step size. In 2D, we have a similar
discretization with $u_{ij}^n = u(t_n, x_i, y_j)$ at the mesh nodes, where  $y_j = (j - 1)
\Delta y$, $j = 1, \dots, N_y$, and $\Delta y$ being the distance between two neighboring
nodes in the $y$-direction.

We use two finite differences schemes to approximately solve the linear advection equation:
a first-order upwind scheme and a second-order Lax-Wendroff scheme. The upwind scheme is
first-order accurate in space and time. In 1D, the discretization of the
linear advection equation is given by
\begin{equation}
    u_i^{n+1} = u_i^n - \frac{a_x \Delta t}{\Delta x}\left(u_i^n - u_{i-1}^n\right),
    \label{eqn:AD1_fd}
\end{equation}
and in 2D by
\begin{equation}
    u_{ij}^{n+1} = u_{ij}^n - \frac{a_x \Delta t}{\Delta x}\left(u_{ij}^n - u_{i-1j}^n\right) - \frac{a_y \Delta t}{\Delta y}\left(u_{ij}^n - u_{ij-1}^n\right).
    \label{eqn:AD2_fd}
\end{equation}

The Lax-Wendroff scheme is second-order accurate in space and time. In 1D, it is given by
\begin{equation}
    u_i^{n+1} = u_i^n - \frac{a_x \Delta t}{2 \Delta x} \left(u_{i+1}^n - u_{i-1}^n\right) + \frac{a_x^2 \Delta t ^ 2}{2 \Delta x^2} \left(u_{i+1}^n - 2 u_i^n + u_{i-1}^n\right)
    \label{eqn:AD1_lw}
\end{equation}
and in 2D by
\begin{equation}
    \begin{split}
        u_{ij}^{n+1} = &\left(1 - \frac{a_x^2 \Delta t^2}{\Delta x^2} - \frac{a_y^2 \Delta t^2}{\Delta y^2}\right) u_{ij}^n + \left(\frac{a_x^2 \Delta t^2}{2 \Delta x^2} - \frac{a_x \Delta t}{2 \Delta x}\right) u_{i+1j}^n \\
        &+ \left(\frac{a_x^2 \Delta t^2}{2 \Delta x^2} + \frac{a_x\Delta t}{2 \Delta x}\right) u_{i-1j}^n + \left(\frac{a_y^2 \Delta t^2}{2 \Delta y^2} - \frac{a_y \Delta t}{2 \Delta y}\right) u_{ij+1}^n \\
        &+ \left(\frac{a_y^2 \Delta t^2}{2 \Delta y^2} + \frac{a_y\Delta t}{2 \Delta y}\right) u_{ij-1}^n\\
        &+ \frac{a_x a_y \Delta t^2}{4 \Delta x \Delta y}\left(u_{i+1j+1}^n - u_{i+1j-1}^n - u_{i-1j+1}^n + u_{i-1j-1}^n\right)
    \end{split}
    \label{eqn:AD2_lw}
\end{equation}

As stated above, we assume periodic boundary conditions in all spatial directions and that
the advection speeds are positive. More details about the employed finite-difference schemes
can be found in, e.g., \cite{LeVeque07}.

\subsection{Reformulating finite difference schemes for FHE}\label{sec:rewrite_fhe}

All FHE operations in the CKKS scheme are vector-based, thus we need to rewrite the finite
difference schemes in Eqns.~\eqref{eqn:AD1_fd}--\eqref{eqn:AD2_lw} in a vector-wise manner.
Furthermore, we need to limit ourselves to the available operations of element-wise
addition, multiplication, and rotation. We begin with the 1D equations and then proceed
to the 2D equations. In the latter case, we also need to develop a strategy to handle
matrix operations with FHE.

In the following, vector/matrix addition and subtraction are
written using the usual symbols $+$ and $-$, respectively. Element-wise vector/matrix
multiplication is indicated by the Hadamard product symbol $\odot$. In the spirit of the
SecureArithmetic.jl library, we do not distinguish in our notation between unencrypted and
encrypted data, since all algorithms work transparently with either.

We represent a generalized rotation operation by $\circshift(\bm{x}, k)$, where a vector
$\bm{x}$ is cyclically shifted by $k$ positions. If $k$ is positive, the shift moves
elements forward; otherwise, it moves them backward. For example,
\begin{equation}
  \circshift\left(\begin{pmatrix} a \\ b \\ c \end{pmatrix}\!\!, 1\right) = \begin{pmatrix} c \\ a \\ b \end{pmatrix}.
\end{equation}
For algorithms in two dimensions, we need to extend the circular shift operation to
matrices. In case of a matrix $\bm{A} \in {\mathbb{R}}^{n \times m}$, $\circshift(\bm{A},
k, l)$ denotes a cyclic shift by $k$ rows and $l$ columns. For example,
\begin{equation}
    \circshift\left(\begin{pmatrix} a & b & c \\ d & e & f \\ g & h & i \end{pmatrix}\!\!, 1, 2\right) = \begin{pmatrix} h & i & g \\ b & c & a \\ e & f & d \end{pmatrix}.
\end{equation}
Depending on the context, the ``$\circshift$'' operator may represent either the rotation of
unencrypted or encrypted vectors and matrices. When we specifically refer to the CKKS
rotation operation, we use the term ``$\rotate$''. For our goal to rewrite the finite
difference schemes in terms of operations supported by FHE, we thus have
addition/subtraction, multiplication, and $\circshift$ at our disposal.

\subsubsection{Rewriting the finite difference methods in 1D}\label{sec:rewrite_1d}

With the definitions above, we can rewrite the 1D first-order upwind finite difference
scheme from Eqn.~\eqref{eqn:AD1_fd} with periodic boundary conditions as
\begin{equation}
    \bm{u}^{n+1} = \bm{u}^n - \frac{a_x \Delta t}{\Delta x}\Bigl(\bm{u}^n - \circshift(\bm{u}^n, 1)\Bigr).
    \label{eqn:AD1_fd_vec_bc}
\end{equation}
Similarly, we can rewrite the 1D Lax-Wendroff scheme from Eqn.~\eqref{eqn:AD1_lw} with
periodic boundary conditions as
\begin{equation}
    \begin{split}
        \bm{u}^{n+1} = \bm{u}^n &- \frac{a_x \Delta t}{2 \Delta x} \Bigl(\circshift(\bm{u}^n, -1) - \circshift(\bm{u}^n, 1)\Bigr)\\
        &+ \frac{a_x^2 \Delta t ^ 2}{2 \Delta x^2} \Bigl(\circshift(\bm{u}^n, -1) - 2 \bm{u}^n + \circshift(\bm{u}^n, 1)\Bigr)
    \end{split}
    \label{eqn:AD1_lw_vec_bc}
\end{equation}

As described in Sec.~\ref{sec:encoding}, the CKKS scheme encodes vectors of real numbers
into plaintexts with a certain batch size/capacity, which is at most half the ring
dimension $N_R$.
If the batch size coincides with the length of the solution vector
$\bm{u}$, we can directly use the CKKS rotation operation for our
$\circshift$ operation. In general, however, we need a way to handle the rotation of vectors
where the length of the user-provided vector is smaller than the capacity. This is
especially relevant for practical applications, since the batch size is always a power
of two.

Therefore, we need to augment the $\circshift$ operation to support the case where the
length of the actual data vector is less than the capacity of the ciphertext. We present an
algorithm for this in Alg.~\ref{alg:circshift_vec}. The algorithm performs two rotations of
the input ciphertext $\bm{x}$, one backward and one forward, and then uses multiplication
with appropriate masking vectors to combine the elements from both rotations. In the
algorithm, we designate the CKKS rotation operation as ``$\rotate$''. Since the OpenFHE
library uses a different sign convention for the rotation index than we use for the
$\circshift$ operation (which follows the sign convention of the corresponding function in
the Julia base library, \texttt{Base.circshift}), we need to adjust the sign of the
shift in the $\rotate$ operation. For the input vectors in
Alg.~\ref{alg:circshift_vec}, we consider the \emph{length} to be equal to the length of the
actual user data, while the \emph{capacity} is equal to the batch size.
\begin{algorithm}[!htbp]
    \caption{Circular shift for secure vectors.}
    \begin{algorithmic}[1]
        \LComment{Input vector $x$, shift index $k$}
        \Function{circshift}{$\bm{x}$, $k$}
            \If{length$(\bm{x}) == $ capacity$(\bm{x})$}
                \State \textbf{return} $\rotate(\bm{x}, -k)$
                \Comment{Use CKKS rotation if data length matches ciphertext capacity}
            \EndIf
            \item[]

            \State $\bm{u} \gets \rotate(\bm{x}, -k)$
            \State $\bm{v} \gets \rotate(\bm{u}, ( k>0$ ? length$(\bm{x})$ : -length$(\bm{x}))-k)$
            \item[]

            \If{$k<0$}
                \Comment{Determine indices for masking vectors}
                \State $f_1 \gets 1$
                \State $l_1 \gets \text{length}(x) + k$
                \State $f_2 \gets 1 + \text{length}(x) + k$
                \State $l_2 \gets \text{length}(x)$
            \Else
                \State $f_1 \gets 1 + k$
                \State $l_1 \gets \text{length}(x)$
                \State $f_2 \gets 1$
                \State $l_2 \gets k$
            \EndIf
            \item[]

            \State $\bm{m}_1 \gets (0, \dots, 0)^\intercal$
            \Comment{Create masking vectors as plaintexts}
            \State $\bm{m}_1[f_1:l_1] \gets 1$
            \State $\bm{m}_2 \gets (0, \dots, 0)^\intercal$
            \State $\bm{m}_2[f_2:l_2] \gets 1$
            \item[]

            \State $\bar{\bm{u}} \gets \bm{u} \odot \bm{m}_1$\Comment{Apply masks by
            multiplying ciphertexts with plaintexts}
            \State $\bar{\bm{v}} \gets \bm{v} \odot \bm{m}_2$
            \item[]

            \State \textbf{return} $\bar{\bm{u}} + \bar{\bm{v}}$
            \Comment{Combine masked ciphertexts for final result}
        \EndFunction
    \end{algorithmic}
    \label{alg:circshift_vec}
\end{algorithm}

With this implementation for the $\circshift$ operation, we now support all operations to
compute the solution at the new time step in Eqns.~\eqref{eqn:AD1_fd_vec_bc} and
\eqref{eqn:AD1_lw_vec_bc}. The $\circshift$ operation consumes one multiplicative level if the
length of the ciphertext is less than its capacity, and zero levels otherwise. Therefore, a
single evaluation of Eqn.~\eqref{eqn:AD1_fd_vec_bc} or Eqn.~\eqref{eqn:AD1_lw_vec_bc} uses
two multiplicative levels if the length of the ciphertext $\bm{u}$ is less than the
capacity, and one level otherwise. Alg.~\ref{alg:circshift_vec} is also the basis for the
$\circshift$ implementation in the SecureArithmetic.jl package.

\subsubsection{Matrix arithmetic with FHE ciphertexts}\label{sec:matrix_fhe}

To reformulate the 2D finite difference schemes in matrix-vector form, we first need to
establish a means of representing data as a matrix in CKKS. Unfortunately, currently no FHE
scheme natively supports matrix arithmetic, thus we need to find a vector-based matrix
representation and define arithmetic operations on it.

The simplest solution would be to put every row of a matrix in a separate ciphertext.
Unfortunately, this approach quickly becomes inordinately inefficient for problems with a
user data length that would otherwise fit into a single ciphertext. Consider, for
example, a 2D domain discretized by $40 \times 40$ nodes. Since the $40 \times 40$ matrix
has only $1600$ elements, it could still be held in a single ciphertext with ring dimension
$N_R = 2^{17}$. However, we can also store the solution matrix as $40$ vectors with a length
of $40$, each in a separate ciphertext. In the latter case, each FHE operation would require
approximately $40$ times longer,
since it would have to be performed for each ciphertext separately. Furthermore, due to
security requirements and other algorithmic constraints, such a storage strategy still
results in ciphertexts with a ring dimension of $2^{17}$. Thus, the size of the matrix in
memory will be around $2.5$ GB for a multiplicative depth of $l_\text{max} = 31$, as opposed
to a single ciphertext with $64$ MB. Obviously, we need to find another vector-based matrix
representation.

\begin{figure}[!htbp]
    \centering
	\includegraphics[width=0.7\linewidth]{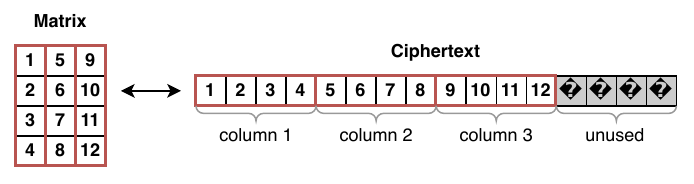} 
    \caption{Storing a $4 \times 3$ matrix (left) in a one-dimensional ciphertext with capacity $16$ (right) by fusing its columns in column-major order. The last four slots of the ciphertext remain unused.}
    \label{fig:matrix_representation}
\end{figure}

Therefore, we store a matrix by fusing its elements column by column and storing it in a
single ciphertext. This is akin to programming languages that store matrix
data linearly in memory in a column-major order. The procedure is visualized in
Fig.~\ref{fig:matrix_representation}. Since the batch size/capacity of CKKS ciphertexts is
always a power of two, in this example the last four slots in the ciphertext are unused and
of indeterminate value.
The good news is that ciphertext addition and multiplication in CKKS are element-wise, and
thus element-wise matrix addition and multiplication are supported out of the box. However,
the $\circshift$ operation does not extend directly to the matrix, since the underlying CKKS
$\rotate$ operation performs only one-dimensional rotations.

We therefore need to extend Alg.~\ref{alg:circshift_vec} to two dimensions. One possible
approach is given in Alg.~\ref{alg:circshift_matrix}. Similar to the 1D algorithm, it works
through a series of rotations and appropriate masking vectors, which are then additively
combined for the final result. With the ``size'' function we obtain the number of rows and
columns of the matrix. The two-dimensional $\circshift$ operation consumes between
zero and two multiplicative levels, depending on whether a mask is required and if the
underlying ciphertext holding the matrix data has a length that is equal or less than its
capacity. An overview of the required number of multiplications is given in
Table~\ref{tab:circshift_2d}. Alg.~\ref{alg:circshift_matrix} is also the basis
for the $\circshift$ implementation for ciphertexts representing matrix data in the
SecureArithmetic.jl package.

\begin{algorithm}[!htbp]
    \caption{Circular shift for secure matrices.}
    \begin{algorithmic}[1]
        \LComment{Input matrix $x$, shift indices $k, l$ for rows and columns.}
        \Function{circshift}{$\bm{x}$, $k$, $l$}
            \State $n, m = \text{size}(\bm{x})$
            \Comment{Get matrix dimensions $n \times m$}
            \If{$k==0$}
            \Comment{If no row shift, directly use vectorial $\circshift$ from Alg.~\ref{alg:circshift_vec}}
                \State \textbf{return} $\circshift(\bm{x}, l n)$
            \EndIf
            \item[]

            \If{$k > 0$}
            \Comment{Determine indices for masking vectors}
                \State $f_1 \gets 1$
                \State $l_1 \gets n - k$
                \State $f_2 \gets 1 + n - k$
                \State $l_2 \gets n$
            \Else
                \State $f_1 \gets 1 - k$
                \State $l_1 \gets n$
                \State $f_2 \gets 1$
                \State $l_2 \gets -k$
            \EndIf
            \item[]

            \State $\bm{m}_{\text{tmp},1} \gets (0,0,\dots,0)^\intercal \in \mathbb{R}^n$
            \Comment{Create masking vectors for one column as plaintexts}
            \State $\bm{m}_{\text{tmp},1}[f_1:l_1] \gets 1$
            \State $\bm{m}_{\text{tmp},2} \gets (0,0,\dots,0)^\intercal \in \mathbb{R}^n$
            \State $\bm{m}_{\text{tmp},2}[f_2:l_2] \gets 1$
            \State $\bm{m}_1 \gets $repeat$(\bm{m}_{\text{tmp},1}, m) \in \mathbb{R}^{nm}$
            \Comment{Repeat the mask for each column}
            \State $\bm{m}_2 \gets $repeat$(\bm{m}_{\text{tmp},2}, m) \in \mathbb{R}^{nm}$
            \item[]

            \State $\bar{\bm{u}} \gets \bm{x} \odot \bm{m}_1$\Comment{Apply masks by
            multiplying ciphertexts with plaintexts}
            \State $\bar{\bm{v}} \gets \bm{x} \odot \bm{m}_2$
            \item[]

            \State $s_1 \gets ln + k$
            \Comment{Compute one-dimensional shift indices}
            \State $s_2 \gets ln + k + (k < 0$ ? $n$ : $-n)$
            \item[]

            \If{$l==0$}
            \Comment{Without column shift, save multiplications by using $\rotate$}
                \State \textbf{return} $\rotate(\bar{\bm{u}}, -s_1) +
                \rotate(\bar{\bm{v}}, -s_2)$
            \Else
                \State \textbf{return} $\circshift(\bar{\bm{u}}, s_1) +
                \circshift(\bar{\bm{v}}, s_2)$
            \EndIf
        \EndFunction
    \end{algorithmic}
    \label{alg:circshift_matrix}
\end{algorithm}

\begin{table}[!htbp]
	\centering
    \begin{subtable}{0.48\textwidth}
        \centering
        \begin{tabular}{ccc}
            \toprule
            $k$  & $l$  & Levels \\
            \midrule
            $= 0$  & $= 0$  & $0$    \\
            $= 0$  & $\neq 0$ & $1$    \\
            $\neq 0$ & $= 0$  & $1$    \\
            $\neq 0$ & $\neq 0$ & $2$    \\
            \bottomrule
        \end{tabular}
        \caption{Length $\neq$ capacity.}
    \end{subtable}
    \begin{subtable}{0.48\textwidth}
        \centering
        \begin{tabular}{ccc}
            \toprule
            $k$  & $l$  & Levels \\
            \midrule
            $= 0$  & $= 0$  & $0$    \\
            $= 0$  & $\neq 0$ & $0$    \\
            $\neq 0$ & $= 0$  & $1$    \\
            $\neq 0$ & $\neq 0$ & $1$    \\
            \bottomrule
        \end{tabular}
        \caption{Length $=$ capacity.}
    \end{subtable}
    \caption{Multiplicative levels required by two-dimensional $\circshift$ for
    rotation by $k$ rows and $l$ columns.}
    \label{tab:circshift_2d}
\end{table}

\subsubsection{Rewriting the finite difference methods in 2D}\label{sec:rewrite_2d}
With the algorithmic building blocks of Sec.~\ref{sec:matrix_fhe}, we can finally rewrite
the 2D upwind scheme from Eqn.~\eqref{eqn:AD2_fd} with periodic boundary conditions in
matrix-vector formulation as
\begin{equation}
    \begin{split}
        \bm{u}^{n+1} = \bm{u}^n &- \frac{a_x \Delta t}{\Delta x} \Bigl(\bm{u}^n - \circshift\bigl(\bm{u}^n, 1, 0\bigr)\Bigr) \\
        &- \frac{a_y \Delta t}{\Delta y} \Bigl(\bm{u}^n - \circshift\bigl(\bm{u}^n, 0, 1\bigr)\Bigr)
    \end{split}
    \label{eqn:AD2_fd_vec}
\end{equation}
Similarly, we can rewrite the 2D Lax-Wendroff scheme from Eqn.~\eqref{eqn:AD2_lw} with periodic boundary conditions as
\begin{equation}
    \begin{split}
        \bm{u}^{n+1} = \Bigl(1 - \frac{a_x^2 \Delta t^2}{\Delta x^2} - \frac{a_y^2 \Delta t^2}{\Delta y^2}\Bigr) \bm{u}^n &+ \Bigl(\frac{a_x^2 \Delta t^2}{2 \Delta x^2} - \frac{a_x \Delta t}{2 \Delta x}\Bigr) \circshift\bigl(\bm{u}^n, -1, 0\bigr) \\
        &+ \Bigl(\frac{a_x^2 \Delta t^2}{2 \Delta x^2} + \frac{a_x\Delta t}{2 \Delta x}\Bigr) \circshift\bigl(\bm{u}^n, 1, 0\bigr) \\
        &+ \Bigl(\frac{a_y^2 \Delta t^2}{2 \Delta y^2} - \frac{a_y \Delta t}{2 \Delta y}\Bigr)  \circshift\bigl(\bm{u}^n, 0, -1\bigr) \\
        &+ \Bigl(\frac{a_y^2 \Delta t^2}{2 \Delta y^2} + \frac{a_y \Delta t}{2 \Delta y}\Bigr) \circshift\bigl(\bm{u}^n, 0, 1\bigr) \\
        &+ \frac{a_x a_y \Delta t^2}{4 \Delta x \Delta y}\Bigl( \circshift\bigl(\bm{u}^n, -1, -1\bigr) - \circshift\bigl(\bm{u}^n, -1, 1\bigr)\\
        &\qquad\qquad\quad \,\,\,-\circshift\bigl(\bm{u}^n, 1, -1\bigr) +  \circshift\bigl(\bm{u}^n, 1, 1\bigr)\Bigr)
    \end{split}
    \label{eqn:AD2_lw_vec}
\end{equation}
$\bm{u}^{n+1}$ is a matrix now, with coefficients $u_{ij}$ at each node location.
Consequently, we also need to use the matrix representation of real data in a ciphertext and
the matrix version of the $\circshift$ operation, which were introduced in the previous
section.

\subsection{Full algorithm for FHE-secured simulations}\label{sec:full_algorithm}

All introduced numerical schemes in Eqns.~\eqref{eqn:AD1_fd_vec_bc}--\eqref{eqn:AD2_lw_vec}
are one-step methods. Alg.~\ref{alg:num_sim} represents the general structure of how all
such numerical methods can be implemented with FHE. Here, it is assumed that the initial
solution $\bm{u}^0$ is already encrypted. The algorithm iterates over the time span $t_0$ to
$t_\text{end}$ with a fixed time step $\Delta t$. In each iteration, we check if during the
next iteration (requiring $l_\text{step}$ multiplicative levels), the ciphertext level would
would fall below one (or two, if iterative bootstrapping is used). If this is the case, the
solution is bootstrapped to increase the ciphertext level again. The actual computation of
the solution at the next time step is then performed in line \ref{alg:line:compute_step}.
The algorithm is terminated when the final time $t_\text{end}$ is reached.

\begin{algorithm}[!htbp]
    \caption{FHE algorithm for secure numerical simulations.}
    \begin{algorithmic}[1]
        \State $\bm{u} \gets \bm{u^0}$
        \Comment{Initialize solution}
        \State $t \gets t_0$
        \Comment{Initialize time}
        \item[]
        
        \While{$t < t_\text{end}$}
            \If{level$(\bm{u}) - l_\text{step} < 1$}
            \Comment{Perform bootstrapping if required}
                \State $\bm{u} \gets $bootstrap$(\bm{u})$
            \EndIf
            \item[]
            
            \State $\bm{u}^{n+1} \gets f(\bm{u})$
            \Comment{Compute solution at next time step based on Eqns.~\eqref{eqn:AD1_fd_vec_bc}--\eqref{eqn:AD2_lw_vec}}\label{alg:line:compute_step}
            \State $\bm{u} \gets \bm{u}^{n+1}$
            \item[]
            
            \State $t \gets t + \Delta t$
            \Comment{Update time}
        \EndWhile
    \end{algorithmic}
    \label{alg:num_sim}
\end{algorithm}

Compared to a classical numerical simulation implementation, the main difference is the
check for the multiplicative depth and the execution of the bootstrapping operation if
required. One should keep in mind that since branching on encrypted ciphertext values is
impossible by design, the condition for exiting a loop cannot be cryptographically secured.
In our implementation we chose to keep the time and step size a plaintext value. To also
hide this information from a potential adversary, one could prescribe a fixed number of
iterations instead, after which the result is returned. This is also the only option if, in
a modified algorithm, the time step were to be computed dynamically from the encrypted
solution data. In addition, we made the decision to further keep some the of the numerical
setup as plaintext, such as the advection speeds or the spatial step size. This is due to
the fact that division is not natively supported by FHE, and thus factors such as $a_x
\Delta t / \Delta x$ need to be computed a priori. However, the actual solution data is
always kept in encrypted form.

The accuracy and performance of cryptographically secure numerical simulations depends
significantly on the performed FHE operations. For Alg.~\ref{alg:num_sim} and the
numerical schemes introduced in the previous sections, we conducted a static analysis of the
number of FHE operations (add, multiply, rotate, and bootstrap) 
that are executed in each iteration. Our findings are
summarized in Table~\ref{tab:num_operations}.
\begin{table}[!htbp]
    \centering
    \begin{tabular}{cccc}
        \multicolumn{4}{c}{length $\neq$ capacity} \\
        \toprule
        \multicolumn{2}{c}{Upwind} & \multicolumn{2}{c}{Lax-Wendroff} \\
        1D  & 2D & 1D  & 2D \\
        \cmidrule(r{0.25em}){1-2} \cmidrule(l{0.25em}){3-4}
        $3$ & $6$ & $4$ & $24$\\
        $3$ & $6$ & $7$ & $38$\\
        $2$ & $4$ & $4$ & $24$\\
        $2/l_\text{usable}$ & $2/l_\text{usable}$ & $2/l_\text{usable}$ & $3/l_\text{usable}$\\
        \bottomrule
    \end{tabular}
    \hspace*{1mm}
    \begin{tabular}{c}
        ~ \\
        \toprule
        CKKS \\
        operation \\
        \midrule
        Add \\
        Multiply\\
        Rotate \\
        Bootstrap\\
        \bottomrule
    \end{tabular}
    \hspace*{1mm}
    \begin{tabular}{cccc}
        \multicolumn{4}{c}{length $=$ capacity} \\
        \toprule
        \multicolumn{2}{c}{Upwind} & \multicolumn{2}{c}{Lax-Wendroff} \\
        1D  & 2D & 1D  & 2D \\
        \cmidrule(r{0.25em}){1-2} \cmidrule(l{0.25em}){3-4}
        $2$ & $5$ & $2$ & $14$\\
        $1$ & $4$ & $3$ & $18$\\
        $1$ & $3$ & $2$ & $14$\\
        $1/l_\text{usable}$ & $2/l_\text{usable}$ & $1/l_\text{usable}$ & $2/l_\text{usable}$\\
        \bottomrule
    \end{tabular}
    \caption{Number of FHE operations required for a  single time step for different
    numerical schemes.}
    \label{tab:num_operations}
\end{table}

As expected, the two-dimensional schemes require more FHE operations than their
one-dimensional counterparts. Similarly, the higher-order Lax-Wendroff scheme requires more
operations than the first-order upwind scheme. The number of operations also depends on the
length of the user data: if it is equal to the ciphertext capacity, some of the
circular shifts reduce to simple CKKS rotations, which necessitate fewer
operations than a full $\circshift$ operation. In addition, using
a data size less than the ciphertext capacity increases the multiplicative depth due to
additional masking operations.

While the higher-order Lax-Wendroff scheme is more accurate in theory, it also uses more FHE
operations per iteration, which increases both runtime and the error introduced by the CKKS
scheme. Therefore, these requirements need to be balanced against each other in practice,
and they likely depend on the specific application.

\section{Numerical results and performance of secure numerical simulations}\label{sec:results}

With all ingredients for FHE-compatible algorithms in place, we can now proceed to
conducting secure numerical simulations. Our investigation aims to compare the accuracy and
performance of the cryptographically secure simulations with their classical counterparts.
For the CKKS setup in OpenFHE, we generally use the same parameters as in
Sec.~\ref{sec:experimental_setup}, but with the available depth after bootstrapping being set
to $l_\text{refresh} = 25$. The ring dimension remains at $N_R = 2^{17}$.
All simulations were
conducted using SecureArithmetic.jl with a single thread unless noted otherwise. The code
for the simulations and the results are available in our reproducibility repository
\cite{reproKholodSchlottkeLakemper24}.

For the numerical simulations presented here, the advection speed of the linear advection
equations is set to $a_x = a_y = 1$. The initial condition is given by $u_0(x) = \sin(2\pi
x)$ in 1D and $u_0(x, y) = \sin(2\pi x) \sin(2\pi y)$ in 2D, and periodic boundary
conditions are imposed. Unless noted otherwise, we use $N = N_x = N_y = 64$ equidistant
nodes in each spatial direction to discretize the computational domain, and a simulation
time span of $t \in [0, 1]$. The time step size $\Delta t$ is determined from the CFL
condition with a CFL number of $0.5$.

In the following, we first conduct convergence tests to verify that the FHE-secured
numerical simulations do not break the convergence properties of the numerical schemes, and
present the numerical results of a full simulation. We then discuss the accuracy and
performance of the secure simulations and compare them to their unencrypted counterparts.
Finally, we briefly explore the potential for speeding up the simulations by leveraging
the multi-threading capabilities of the OpenFHE library.

\subsection{Convergence test}\label{sec:convergence_test}

We begin by performing a convergence study to verify the implementation of the numerical
methods and to make sure that the convergence properties of the numerical schemes are not
affected by the FHE encryption. We use the simulation setup as described, but vary the
number of nodes in each direction from $N = 32$ to $N = 256$ and reduce the simulation time
span to $t \in [0, 0.5]$.
To speed up the convergence tests, we ran the simulations in parallel with eight threads
(for more details on parallelizing OpenFHE with OpenMP, see Sec.~\ref{sec:parallelization}
below).

We use the discrete $L^2$ norm as a measure of the error between the numerical solution $u$
and the exact solution $u_\text{exact}$, where $u_{\text{exact},i}^n = u_0(x_i - a_x t_n)$
in 1D and $u_{\text{exact},ij}^n = u_0(x_i - a_x t_n, y_j - a_y t_n)$ in 2D. The $L^2$ error
for $N$ nodes per direction is then calculated by
\begin{equation}
    \text{1D: } \, e_N = \sqrt{\frac{1}{N}\sum_{i=1}^{N}(u_{i}^{n}-u_{\text{exact}, i}^{n})^2} \qquad
    \text{2D: } \, e_N = \sqrt{\frac{1}{N^2}\sum_{i,j=1}^{N}(u_{ij}^{n}-u_{\text{exact}, ij}^{n})^2}
    \label{eqn:l2_error}
\end{equation}
From the $L^2$ error, we then compute the experimental order of convergence (EOC)
at each step by
\begin{equation}
    \text{EOC} = \log_2(e_N/e_{2N})
    \label{EOC}
\end{equation}

Table~\ref{tab:convergence_test} contains the results of the convergence study. The EOC
values are close to the expected values of one for the
first-order upwind scheme and two for the second-order Lax-Wendroff scheme. This confirms
the correctness of the implementation and shows
that the FHE encryption does not break the convergence properties of the numerical schemes.
\begin{table}[!htbp]
	\centering
    \begin{adjustbox}{width=\textwidth}
        \begin{tabular}{rrr}
            \multicolumn{3}{c}{1D, upwind} \\
            \toprule
            $N$ & $L^2$ error & EOC \\
            \midrule
            32 & 1.01e-01 & - \\
            64 & 5.25e-02 & 0.95 \\
            128 & 2.67e-02 & 0.97 \\
            256 & 1.35e-02 & 0.99 \\
            \bottomrule
        \end{tabular}
        \hspace{1mm}
        \begin{tabular}{rrr}
            \multicolumn{3}{c}{1D, Lax-Wendroff} \\
            \toprule
            $N$ & $L^2$ error & EOC \\
            \midrule
            32 & 1.07e-02 & - \\
            64 & 2.67e-03 & 2.00 \\
            128 & 6.69e-04 & 2.00 \\
            256 & 1.67e-04 & 2.00 \\
            \bottomrule
        \end{tabular}
        \hspace{1mm}
        \begin{tabular}{rrr}
            \multicolumn{3}{c}{2D, upwind} \\
            \toprule
                $N$ & $L^2$ error & EOC \\
            \midrule
            32 & 1.88e-01 & - \\
            64 & 1.07e-01 & 0.82 \\
            128 & 5.69e-02 & 0.90 \\
            256 & 2.94e-02 & 0.95 \\
            \bottomrule
        \end{tabular}
        \hspace{1mm}
        \begin{tabular}{rrr}
            \multicolumn{3}{c}{2D, Lax-Wendroff} \\
            \toprule
                $N$ & $L^2$ error & EOC \\
            \midrule
            32 & 1.07e-02 & - \\
            64 & 2.68e-03 & 2.00 \\
            128 & 6.69e-04 & 2.00 \\
            256 & 1.67e-04 & 2.00 \\
            \bottomrule
        \end{tabular}
    \end{adjustbox}
	\caption{
		Convergence test results for the secure simulation of the 1D/2D
		linear scalar advection equation with the first-order upwind scheme and the
		second-order Lax-Wendroff scheme.
	}
  	\label{tab:convergence_test}
\end{table}

\subsection{Results for a full simulation}

\begin{figure}[!htbp]
    \centering
    \begin{subfigure}[b]{0.48\textwidth}
        \centering
        \includegraphics[width=\textwidth]{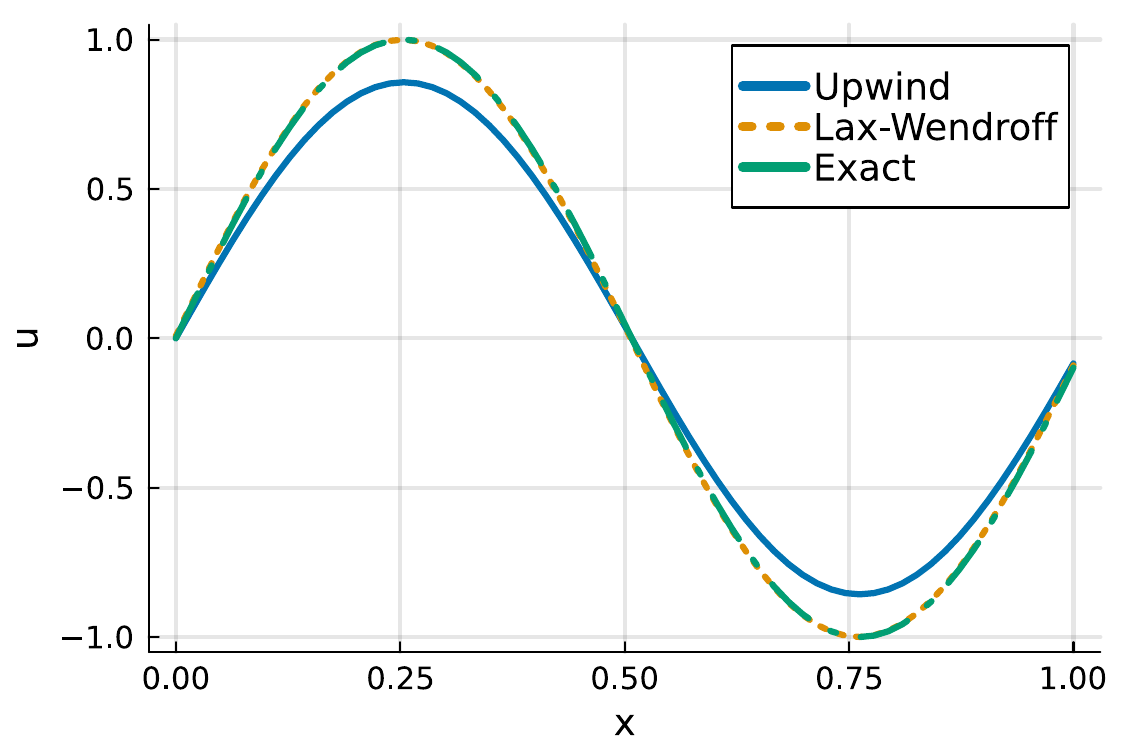}
    \end{subfigure}
    \caption{Secure simulation results of the linear scalar advection equation with the first-order upwind and the second-order Lax-Wendroff schemes in 1D at $t = 1$.}
    \label{fig:vis1d}
\end{figure}

Fig.~\ref{fig:vis1d} shows the results for a full simulation of the linear scalar advection
equation
in 1D at $t = 1$, i.e., after one full period of this periodic problem setup. The results of
the upwind scheme and the Lax-Wendroff scheme are compared to the exact solution. As
expected, the first-order upwind scheme is much more dissipative than the second-order
Lax-Wendroff scheme, with the latter being visually indistinguishable from the exact
solution.

\begin{figure}[!htbp]
    \centering
    \begin{subfigure}[b]{0.48\textwidth}
        \centering
        \includegraphics[width=1.0\textwidth,trim=4cm 1.0cm 3.0cm 1.0cm, clip]{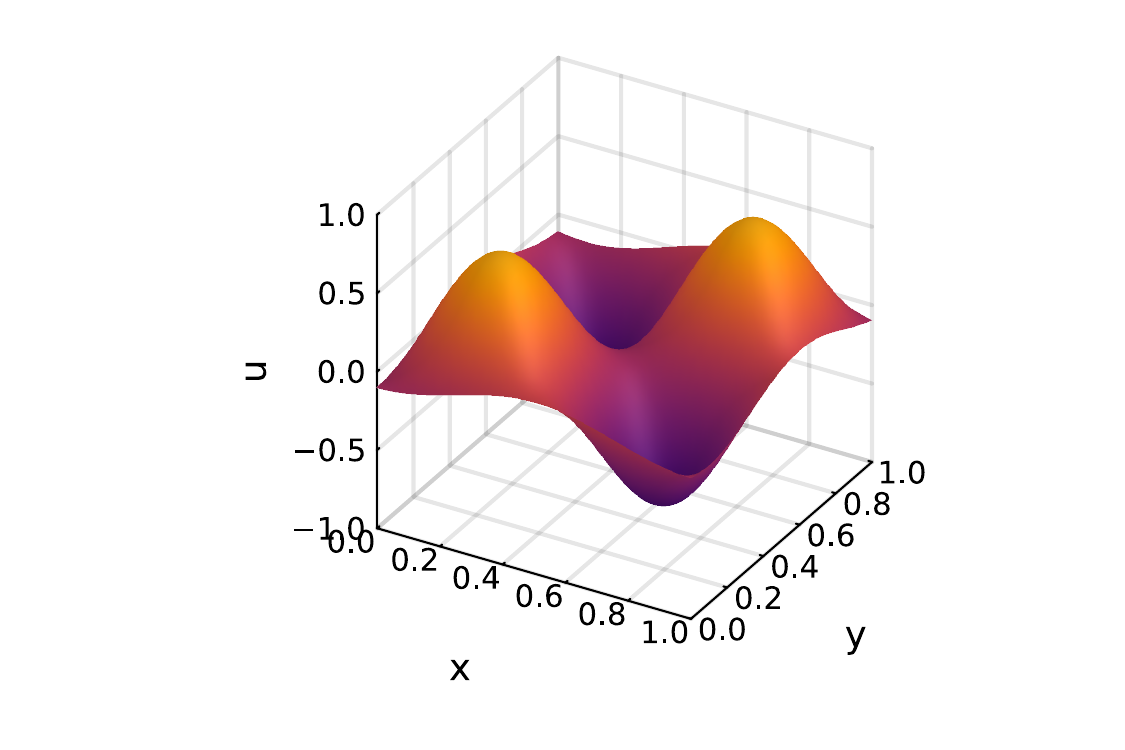}
        \caption{First-order upwind scheme.}
        \label{fig:vis2d_upwind}
    \end{subfigure}%
    \hfill%
    \begin{subfigure}[b]{0.48\textwidth}
        \centering
        \includegraphics[width=1.0\textwidth,trim=4cm 1.0cm 3.0cm 1.0cm, clip]{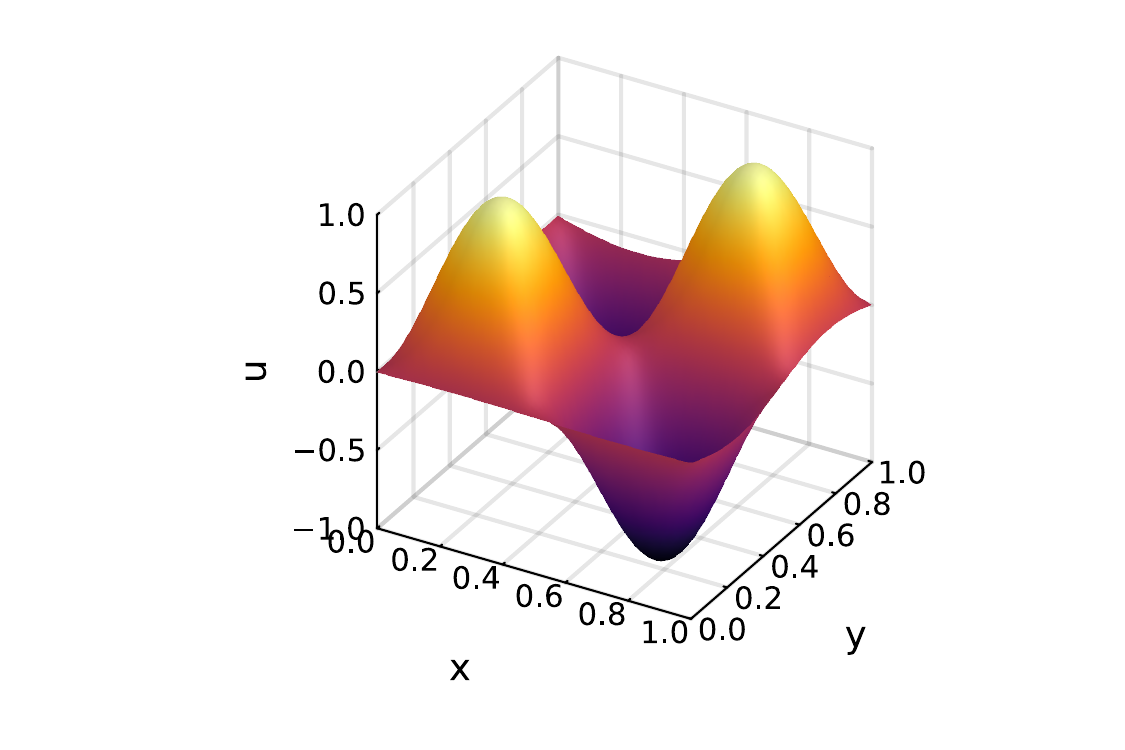}
        \caption{Second-order Lax-Wendroff scheme.}
        \label{fig:vis2d_lw}
    \end{subfigure}
    \caption{Secure simulation results of the linear scalar advection equation in 2D at $t = 1$.}
    \label{fig:vis2d}
\end{figure}

In Fig.~\ref{fig:vis2d}, we present the simulation results for the linear scalar advection
in 2D, also at $t = 1$. Again, the first-order upwind scheme displays significant
dissipative effects. We omitted the exact solution, since it is again indiscernible from the
Lax-Wendroff scheme.

\subsection{Accuracy and performance of encrypted vs.\ unencrypted simulations}\label{sec:accuracy-performance}

\begin{figure}[!htbp]
    \begin{subfigure}{0.49\textwidth}
	  \includegraphics[width=\textwidth]{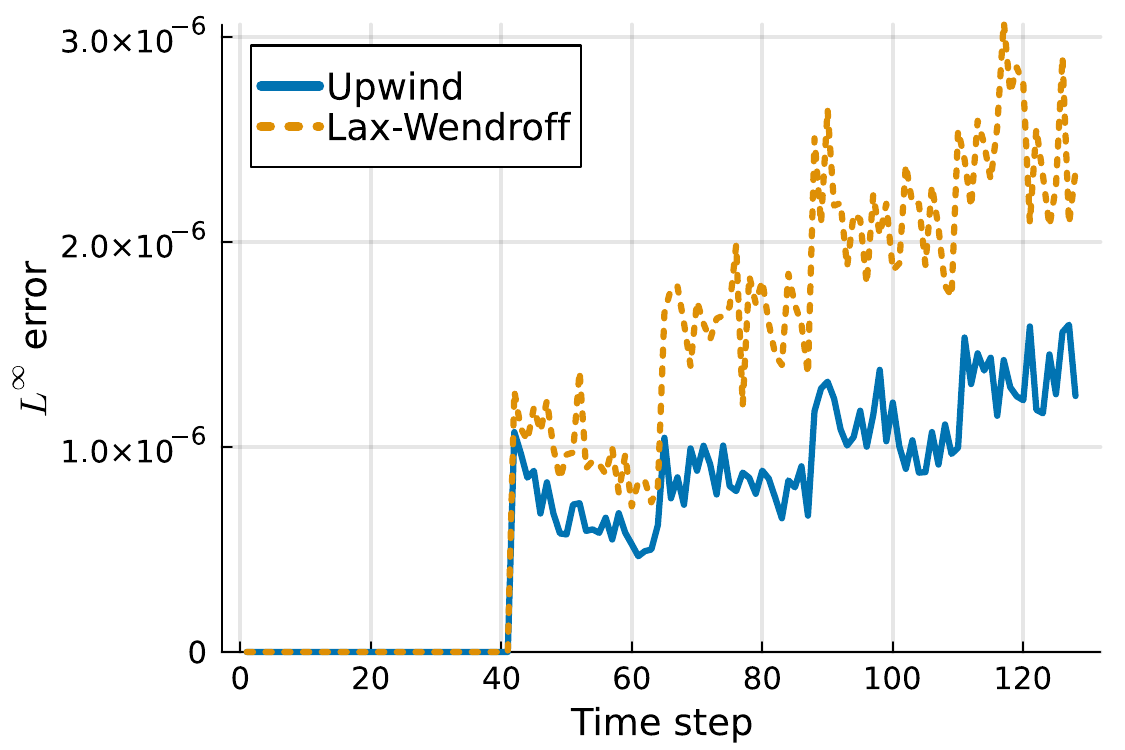}
        \caption{1D.}
    \end{subfigure}
    \hfill
    \begin{subfigure}{0.49\textwidth}
        \includegraphics[width=\textwidth]{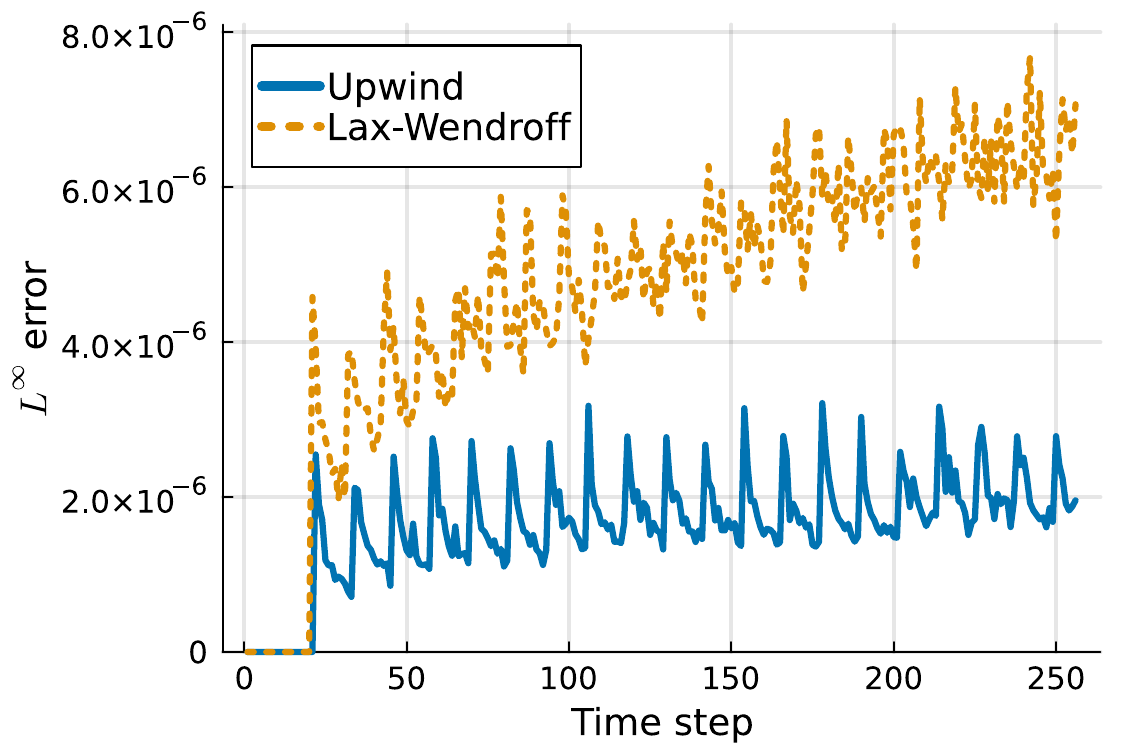}
        \caption{2D.}
    \end{subfigure}
    \caption{Errors introduced by the CKKS scheme for the numerical simulations with different discretization methods.}
    \label{fig:error}
\end{figure}

With the convergence properties of the numerical solvers verified, we proceed with
comparing the accuracy of the results and the performance of the secure algorithms with
their unencrypted counterparts.
In Fig.~\ref{fig:error} we analyze the accuracy of the cryptographically secure numerical
simulations compared to the unencrypted simulations, using the same setup as before. Unlike the convergence test,
where we used the exact solution as a reference, we now compare the results of the secure
simulations to the results of the unencrypted simulations. This allows us to separate the
error of the numerical discretization and instead only consider the error introduced by the
FHE encryption, which we compute as the $L^\infty$ norm of the difference between the
secure and unencrypted solutions at each time step. In both the 1D and 2D cases, we can see
that the error for the Lax-Wendroff scheme is higher than for the upwind scheme.
This is expected, since the Lax-Wendroff scheme requires more FHE operations, which in turn
increases the error introduced by the CKKS scheme. As we saw above, however, overall the
Lax-Wendroff scheme is still much more accurate than the upwind scheme, even with the
additional error introduced by the encrypted operations. Furthermore, the relative
$L^\infty$ error at $\mathcal{O}(10^{-6})$ is still several orders of magnitude lower then
the absolute numerical error of the unencrypted simulations at $\mathcal{O}(10^{-3})$.
Yet, it is also clear that with FHE, less than single-precision floating point accuracy
is achieved, which may be a limiting factor for some applications.

It is also very interesting to see the large effect the bootstrapping operation has on the
accuracy. Until the first bootstrapping operation at around 40 time steps in the 1D
simulations and at around 25 time steps in the 2D simulations, the error is very low at
$\mathcal{O}(10^{-13})$. After the bootstrapping operation, the error jumps by several
orders of magnitude to $\mathcal{O}(10^{-6})$. We also observe that the error seems to
gradually decrease again after the bootstrapping operation, which is especially notable for
the 2D upwind solution. The observation that the precision of CKKS with bootstrapping improves over time is a relatively common phenomenon in privacy-preserving machine learning applications of CKKS. For instance, Han $\textit{et al.}$ show that the CKKS error in encrypted logistic regression training is higher in initial iterations, but then the encrypted solution gradually converges to the plaintext result due to the convergence of the underlying (logistic regression training) procedure (see Section 5.2 of~\cite{cryptoeprint:2018/662}). We expect similar behavior in converging finite difference schemes, and consider the observed precision improvement after bootstrapping as an indicator of this more general phenomenon.

Next, we look at the runtime performance of the secure numerical simulations. In
Table~\ref{tab:runtime_per_step} we show the runtime per step for the different schemes,
which in case of the secure simulations was measured between the first and the second
bootstrapping operation, including bootstrapping, and then averaged over the number of time
steps. From the table, we can clearly see that the secure simulations are significantly
slower than the unencrypted simulations, by approximately six orders of magnitude.
Considering the runtime measurements of the individual CKKS operations in
Sec.~\ref{sec:accuracy-performance}, this is not surprising. Furthermore, the 2D simulations
are more expensive than the 1D simulations, and the Lax-Wendroff scheme has a higher runtime
than the upwind scheme, which also matches expectations. It is interesting to note that
the runtime of the encrypted simulations only differs by a factor of two to four between the
1D and 2D simulations, while the unencrypted 2D simulations are around ten times more
expensive than their 1D counterparts.
Given that the FHE operation count in the 2D simulations is considerably higher than in 1D
(Table~\ref{tab:num_operations}), the relatively modest increase in runtime suggests that
bootstrapping dominates the overall performance cost.
This is further supported by the
measured bootstrapping time: for the Lax-Wendroff scheme, bootstrapping accounts for about
54\% of the total runtime in 1D and 36\% in 2D.
\begin{table}[!htbp]
	\centering
    \begin{subtable}{0.48\textwidth}
        \begin{tabular}{llrr}
            \toprule
            & & OpenFHE & Unencrypted \\
            \midrule
            \multirow{2}{*}{1D} & Upwind & 4.67 s & 6.85e-7 s \\
            & Lax-Wendroff & 5.28 s & 1.21e-6 s \\
            \cmidrule{1-2}
            \multirow{2}{*}{2D} & Upwind & 10.5 s & 2.41e-5 s \\
            & Lax-Wendroff & 20.8 s & 6.92e-5 s \\
            \bottomrule
        \end{tabular}
        \caption{Runtime per time step.}
        \label{tab:runtime_per_step}
    \end{subtable}
    \begin{subtable}{0.48\textwidth}
        \begin{tabular}{llrr}
            \toprule
            & & OpenFHE & Unencrypted \\
            \midrule
            \multirow{2}{*}{1D} & Upwind & 93.3 s & 4.44e-6 s \\
            & Lax-Wendroff & 101.0 s & 5.89e-6 s \\
            \cmidrule{1-2}
            \multirow{2}{*}{2D} & Upwind & 170.1 s & 9.73e-4 s \\
            & Lax-Wendroff & 220.8 s & 1.09e-3 s \\
            \bottomrule
        \end{tabular}
        \caption{Time for initialization.}
        \label{tab:configuration_time}
    \end{subtable}
	\caption{ Runtime per time step (left) and initial configuration time (right) for the 1D
		and 2D simulations with the secure and insecure backends.}
	\label{tab:exec_time}
\end{table}

A brief analysis of the time required for the initialization of the CKKS scheme is given in
Table~\ref{tab:configuration_time}. Not surprisingly, the initial setup of the FHE
operations is orders of magnitude slower than the initialization of the unencrypted data
structures. Even though ciphertexts with the same ring dimensions are used, the 2D times are
much higher than for the 1D simulations, owing to the fact that they have a significantly
larger batch size/capacity. In addition, the 2D $\circshift$ operation requires a larger
number of index shifts, which need to be set up during the initialization phase. Similarly,
the Lax-Wendroff scheme requires rotations in more directions than the upwind scheme, which
also increases the initialization time.

\begin{figure}[!htbp]
    \begin{subfigure}{0.48\textwidth}
	  \includegraphics[width=\textwidth]{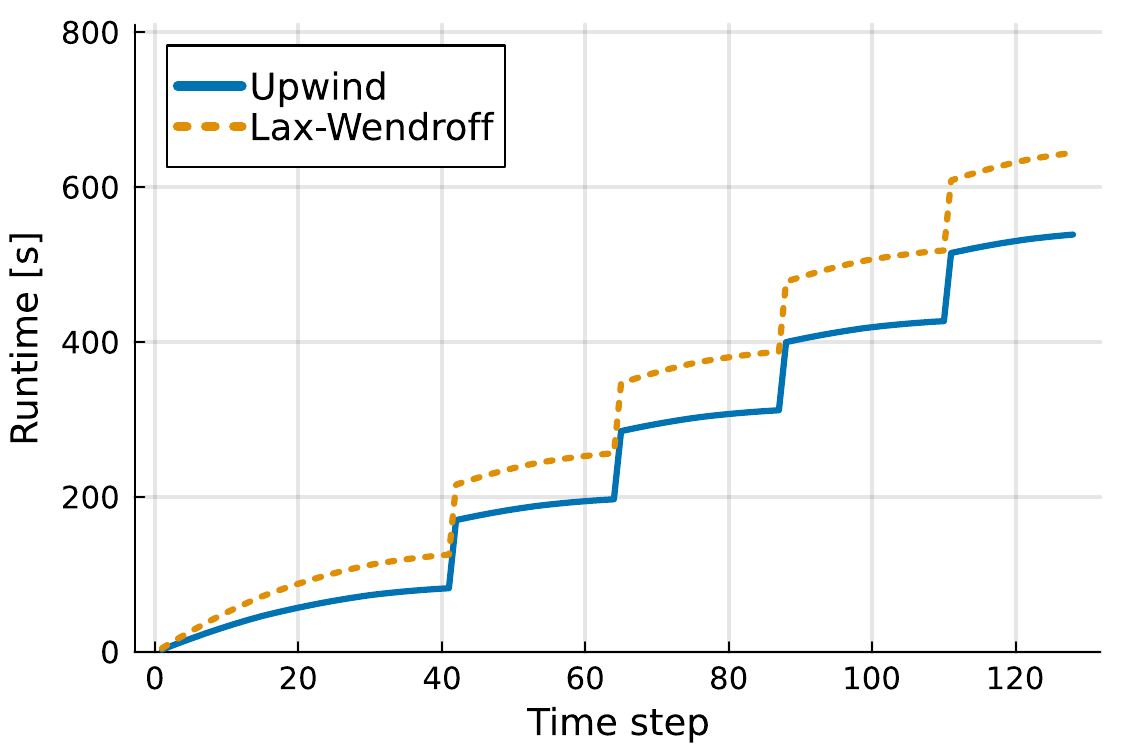}
      \caption{1D.}
      \label{fig:exec_time_1d}
    \end{subfigure}
    \hfill
    \begin{subfigure}{0.48\textwidth}
	  \includegraphics[width=\textwidth]{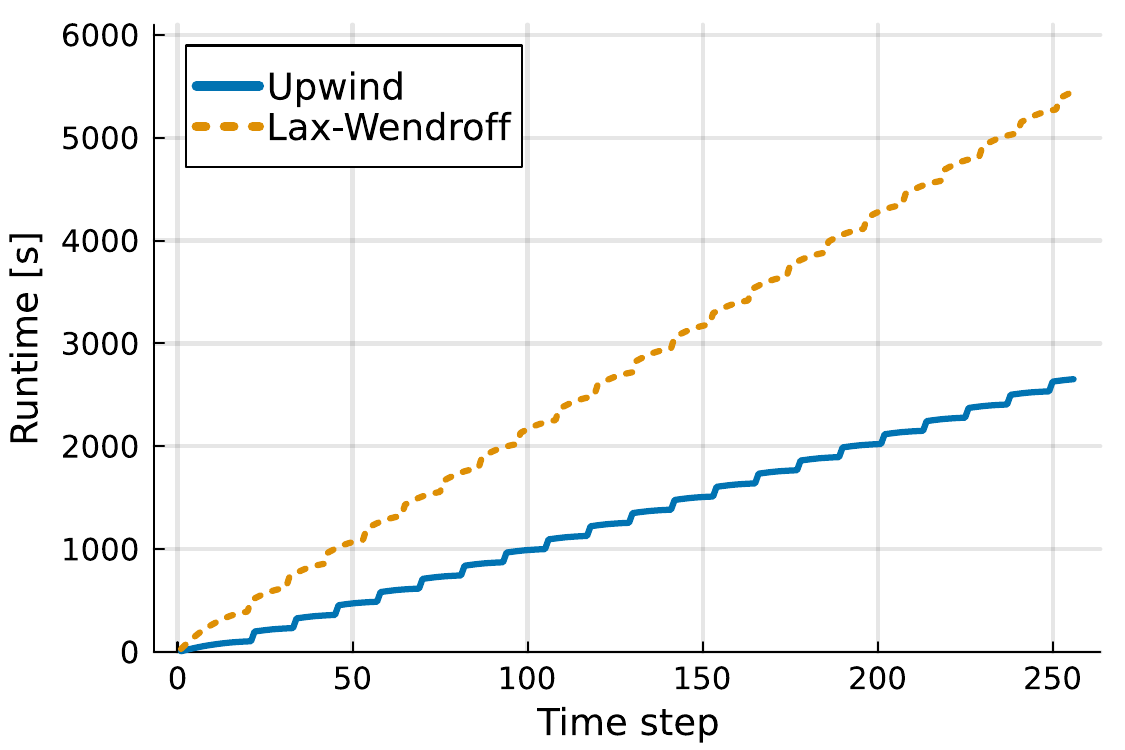}
       \caption{2D.}
       \label{fig:exec_time_2d}
    \end{subfigure}
    \caption{Execution times of secure numerical simulations.}
    \label{fig:exec_time}
\end{figure}

Next to the average runtime per step, it is instructive to look at the runtime
over the number of time steps (Fig.~\ref{fig:exec_time}). Especially for the 1D
simulations in Fig.~\ref{fig:exec_time_1d}, we can see that the run time per step
decreases with the number of steps, until there is a sudden and significant increase in
runtime.
These jumps can be attributed to the bootstrapping operation,
matching the first bootstrapping after around 40 steps and then
every approximately 25 steps thereafter. This also lines up with the jumps in the error at
the same time steps we have seen in Fig.~\ref{fig:error}. The fact that the runtime per step
decreases with each step up to
the next bootstrapping operation also matches the performance characteristics we saw earlier
in Sec.~\ref{sec:accuracy_performance}, where many FHE operations became slower with
increasing multiplicative depth.
This is due to the ciphertext modulus being reduced after each multiplication, which
in turn reduces the number of computational steps required for each subsequent FHE
operation.
The same behavior can be observed here, where the remaining
multiplicative depth decreases with each time step until it is refreshed by
bootstrapping. A similar trend can be observed in Fig.~\ref{fig:exec_time_2d} for the 2D
simulations.

\begin{figure}[!htbp]
    \centering
	\begin{subfigure}{0.48\textwidth}
        \centering
		\includegraphics[width=\textwidth]{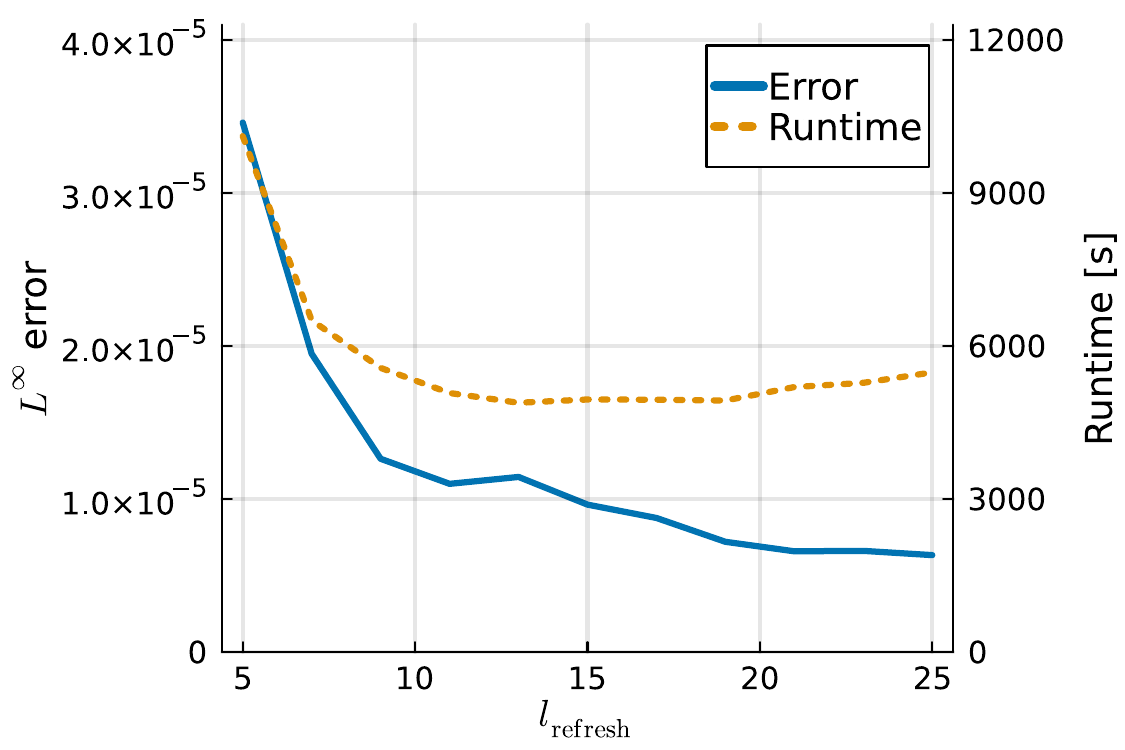}
	\end{subfigure}
	\caption{Error and execution time for secure simulations with different values of
	$l_\text{refresh}$ (multiplicative depths after bootstrapping). The simulations of
	the 2D linear advection equation were performed with the Lax-Wendroff scheme.}
	\label{fig:time_error}
\end{figure}

To highlight their interdependence,
in Fig.~\ref{fig:time_error} we analyze the error and runtime of the secure
simulations with the Lax-Wendroff scheme in 2D for different multiplicative depths after
bootstrapping $l_\text{refresh}$. The error becomes smaller with larger $l_\text{refresh}$,
aligning with our previous observation that the bootstrapping incurs a high approximation
error.
The overall runtime first decreases with increasing multiplicative depth, which is expected
since the bootstrapping operation occurs less frequently.
However, there seems to exist
an optimum value, since after approximately $l_\text{refresh} = 17$ the runtime increases
again. This confirms our earlier discussion that there exists a trade-off between the
computational cost of the FHE operations, which increases with larger values for
$l_\text{refresh}$, and the frequency of the bootstrapping operation. The optimal
multiplicative depth thus depends on the specific application and the used numerical scheme.

\begin{figure}[!htbp]
    \begin{subfigure}{0.48\textwidth}
	  \includegraphics[width=\textwidth]{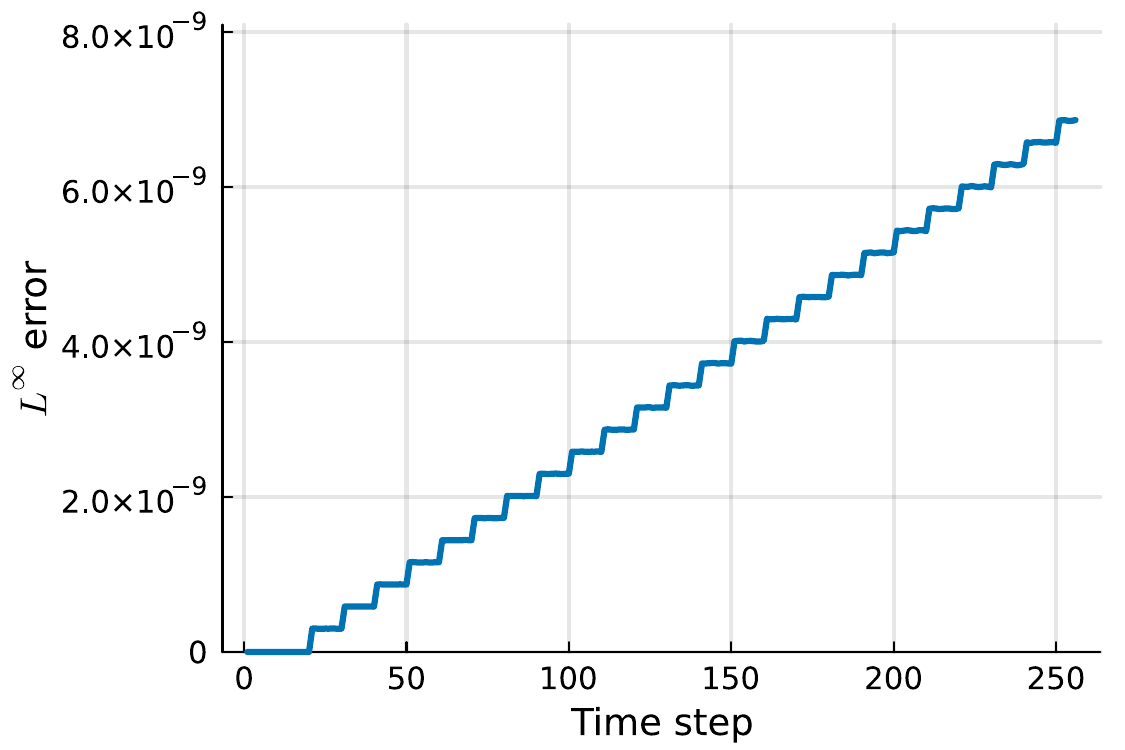}
    \end{subfigure}
    \hfill
    \begin{subfigure}{0.48\textwidth}
	  \includegraphics[width=\textwidth]{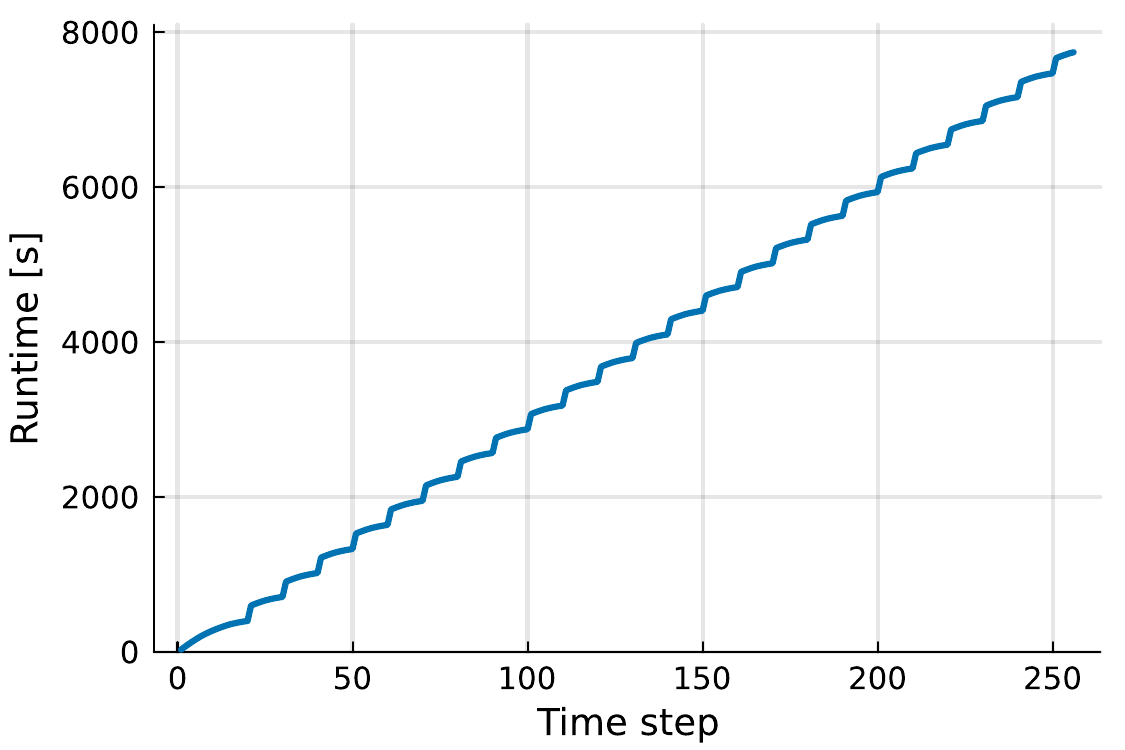}
    \end{subfigure}
    \caption{Error (left) and execution time (right) for secure simulations with iterative
	bootstrapping. The simulations of the 2D linear advection equation were performed with
	the Lax-Wendroff scheme.}
    \label{fig:iterative-bootstrapping-lw2d}
\end{figure}

Finally, we examine the effect of using the iterative bootstrapping technique with two
iterations as discussed in Sec.~\ref{sec:accuracy-performance-bootstrapping}.
Compared to the previous results with standard bootstrapping in Fig.~\ref{fig:error},
Fig.~\ref{fig:iterative-bootstrapping-lw2d} shows that the
error is reduced by three orders of magnitude, from approximately
$\mathcal{O}(10^{-6})$ to $\mathcal{O}(10^{-9})$. This matches our previous findings in
Fig.~\ref{fig:bootstrapping}, reemphasizing the dominant role of the bootstrapping
operation on the overall error caused by the CKKS scheme.
At the same time, the runtime per step
rises from $20.8$~s to $30.4$~s. The increase in runtime is mainly due to doubling
the number of bootstrapping operations. Furthermore, for iterative bootstrapping, a larger
number of available levels $l = 2$ is required (compared to $l = 1$
before, see also Sec.~\ref{sec:accuracy-performance-bootstrapping}). Nevertheless, it is
clear that iterative bootstrapping can be a viable option in scenarios where the error
introduced by the CKKS scheme is a limiting factor and the computational cost is not
prohibitive. Interestingly, the previously observed self-healing capabilities of the finite
difference scheme seem to disappear with iterative bootstrapping, as the error no longer
decreases between two bootstrapping operations.

Overall, our findings demonstrate that secure numerical simulations with the CKKS scheme are
feasible. Results remain accurate as long as bootstrapping is not required, which introduces
noticeable numerical errors---though single-precision accuracy is still achievable. The
results also highlight the considerable performance overhead introduced by CKKS---primarily
from bootstrapping.
For practical applications, it thus seems advisable to primarily
focus on optimizing the multiplicative depth and the bootstrapping configuration to minimize
the introduced error and to maximize the computational performance of secure simulations.
In scenarios where runtime is not the primary limitation, iterative
bootstrapping may be used to improve the accuracy of FHE-based simulations.

\subsection{Parallel computations}\label{sec:parallelization}

As evident from the previous sections, FHE operations incur a considerable performance overhead. Therefore, it makes sense to evaluate potential ways to speed up the secure
computations. One possible approach is to parallelize the FHE operations.

The OpenFHE library internally uses the OpenMP library for multithreading computations,
which is also available in Julia. To study the efficiency of CKKS parallelization, we consider the most expensive numerical solver from this work, i.e., the Lax-Wendroff scheme
for the linear advection equation in 2D, and execute it for various numbers of OpenMP
threads. This can be controlled by setting the environment variable
\texttt{OMP\_NUM\_THREADS} to the appropriate number of threads before starting the Julia runtime.

\begin{figure}[!htbp]
    \begin{subfigure}{0.49\textwidth}
	  \includegraphics[width=\textwidth]{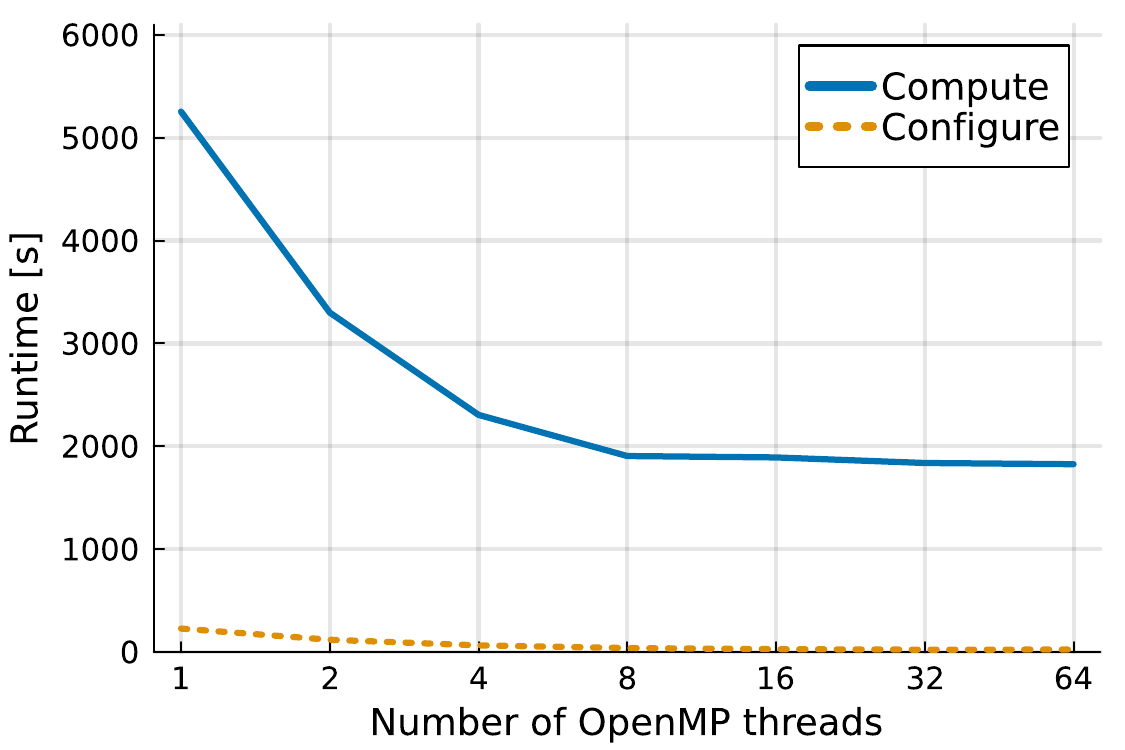}
      \caption{Execution time for different numbers of OpenMP threads.}
      \label{fig:parallel_LW_2D_runtime}
    \end{subfigure}
    \hfill
    \begin{subfigure}{0.49\textwidth}
      \includegraphics[width=\textwidth]{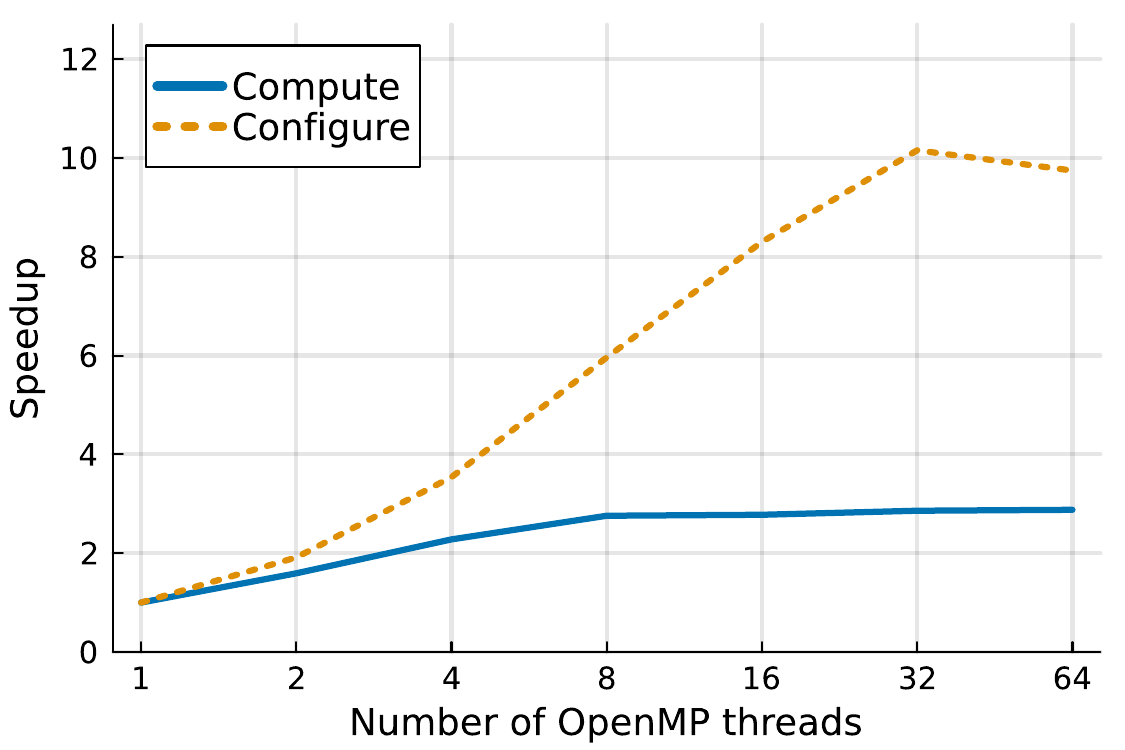}
      \caption{Speedup for different numbers of OpenMP threads.}
      \label{fig:parallel_LW_2D_speedup}
    \end{subfigure}
    \caption{Parallel execution with multiple OpenMP threads of the Lax-Wendroff scheme for the linear scalar advection equation in 2D. Runtime and speedup are analyzed separately for the configuration and the computation phase.}
    \label{fig:parallel_LW_2D}
\end{figure}

Fig.~\ref{fig:parallel_LW_2D} shows the runtimes for the parallel execution of the secure
simulations. For the used CKKS parameters and hardware setup, we could achieve
at most a three-fold speedup in computation time with eight threads, and a ten-fold speedup
with $32$ threads in configuration time. In either case, using more threads than this did
not lead to
any significant improvements in runtime.

Given that using multiple threads can reduce solution times and that the OpenFHE library
offers basic multithreading capabilities, it is advisable to employ 2 to 8 threads to speed
up CKKS computations whenever possible.
However, the observed speedup remains fall short of ideal scaling, indicating that the current
multithreading implementation has practical limitations. Since efficient
parallelization is crucial for the viability of FHE-based simulations---especially in more demanding applications---this remains an important subject for future investigation.
Nevertheless, these results provide a basis for exploring more advanced simulation scenarios, outlined in the next section.

\section{Extension to other mathematical models and numerical methods}\label{sec:extension}
Building on the prototype simulations and findings from Sec.~\ref{sec:results}, we discuss
challenges and possible strategies for extending secure numerical simulations to more
complex models and schemes.
As noted before, the primary impediments to the practical use of the CKKS
scheme for such scenarios are the large computational overhead and, to a lesser extent, the
precision loss. In addition, three main challenges arise when considering to extend the
simulation prototypes shown in this paper to other PDEs and advanced numerical
schemes: a single built-in data structure, limited function evaluation capabilities, and
constraints in control flow logic when it comes to conditional branching. In
the following, we discuss these issues and potential solutions to them, to give other
researchers a starting point for their own investigations into secure numerical simulations.
We then outline general strategies to improve performance and conclude with a brief summary
of simulation types that are (currently) out of scope.

\subsection{Data structures}\label{sec:data_structures}
A vector of real numbers is the only encrypted data structure available in the CKKS scheme,
as described in Sec.~\ref{sec:ckks}. More complex data structures, such as tensors or
matrices, are not natively supported. In Sec.~\ref{sec:matrix_fhe}, we have introduced one
potential approach to a matrix representation for use with explicit, mesh-based numerical
simulations. In other algorithms that involve matrix-matrix or matrix-vector
multiplications, specialized encoding strategies for the matrix lead to more efficient
algorithms in terms of the multiplicative depth \cite{HaleviShoup14, jiang2018matrix, BlattGusevEtAl20}.

In general, the concrete choice of data representation depends on the target application and
the used numerical method. Due to the significant computational overhead of the CKKS scheme,
it is necessary to make full use of its batch processing capabilities, i.e., the ability to
do element-wise addition and multiplication on the entire data vector. Therefore,
discretization schemes that allow for efficient vectorization are preferable when
considering other numerical methods for secure simulations. The choice of PDE is largely
unaffected by the data structure, as long as a single model with a constant number of state
variables is used throughout the computational domain.
However, algorithms that require fine-grained access to the data---such as
irregular sparse matrix operations---are not directly compatible with FHE.

\subsection{Function evaluation}
As discussed before, the CKKS scheme supports only a limited set of arithmetic operations,
i.e., element-wise addition and multiplication of the ciphertexts, as well as
rotation. This means that more complex arithmetic operations, such as division,
exponentiation, or trigonometric functions, are not directly available but must be
approximated by a series of simpler operations or using a polynomial interpolation. For example, the division of two ciphertexts
can be implemented with the Newton-Raphson method \cite{CetinDorozEtAl15} or Goldschmidt's
algorithm \cite{MoonOmarovEtAl24}, the latter of which is also used to implement floating
point division in some CPU families \cite{oberman1999floating}. Exponentiation or
trigonometric functions may be approximated by a Taylor series expansion
\cite{prantl2024debello}.

Since evaluating a polynomial expression is straightforward with the built-in arithmetic
operations, any function that is reasonably well approximated by a polynomial can be
computed with sufficient accuracy. For many other functions, an iterative approach is
required. Here, the main difficulty is that the number of iterations must be fixed a priori,
since branching based on encrypted data is not possible (refer to the next section for
details). Without fundamental limitations to computing branch-free functions with FHE, the
primary challenge is therefore to find an efficient approximation that converges robustly.
Selecting the best approach is problem-dependent and likely requires some experimentation to
find a good balance between accuracy and performance. When considering other mathematical
models or numerical discretizations, it is thus advisable to choose options that minimize
the multiplicative depth of the FHE operations to avoid costly bootstrapping operations.

For many discontinuous functions, e.g., sign/comparison, the use of CKKS is challenging as
the polynomial approximation in the proximity of the singularity point(s) requires extremely
high-degree polynomials, i.e., very large multiplicative depth and runtime. Instead, other
FHE schemes are often used, such as DM \cite{DM15} or CGGI \cite{ChillottiGamaEtAl19}. DM
and CGGI support the evaluation of arbitrary functions via lookup tables using a technique
called functional or programmable bootstrapping~\cite{LMP22}. In these scenarios, a
\emph{scheme switching} procedure is employed \cite{boura2018chimera, lu2020pegasus}, where
the encryption of the data is temporarily changed from CKKS to DM/CGGI. This allows for the
evaluation of the discontinuous function, after which the data is converted back to a CKKS
representation. Such scheme switching capability is supported in
OpenFHE~\cite{OpenFHE}. The main drawback of scheme switching is the significant
computational overhead, especially when high precision is required
\cite{eldefrawy2023hardness, badawi2023demystifying, kim2024grass}. A promising new method
for evaluating discontinuous functions using functional CKKS bootstrapping is proposed and
implemented in \cite{alexandru2024functionalbootstrapping}. The method enables arbitrary
function evaluation via lookup tables and achieves a throughput that is orders of magnitude
higher than DM or CGGI.

\subsection{Data-dependent control flow}\label{sec:control_flow}
By their very nature, many of the control flow directives commonly used in regular programs
cannot be employed with FHE methods. More specifically, FHE precludes the use of conditional
statements or branches---including loops with a termination condition---that depend on
encrypted data. This is due to the fact that the encryption of the data prevents the
processor from knowing the actual values of the data. If this were not the case, it
would be possible for a malicious actor to infer information about the encrypted data,
e.g., by using a bisection algorithm to compare the secret ciphertext with known values.

As noted in \cite{chialva2018conditionals}, this restriction in control flow logic presents
a fundamental challenge when converting regular programs to FHE algorithms. Since
branching on an encrypted value is impossible by design, an FHE algorithm encountering a
conditional branch must therefore first compute the outcome for each possible branch. The
final result is then constructed by combining the results from all branches and applying an
appropriate mask that selects the correct result. For example,
Alg.~\ref{alg:branching_classic} shows a simple branching statement in a regular program.
Based on a function \texttt{condition()} that returns either $1$ or $0$, the value of $r$ is
computed by calling either \texttt{foo()} or \texttt{bar()}. The equivalent FHE algorithm in
Alg.~\ref{alg:branching_fhe} first computes the value of the condition, evaluates both
branches, and then combines the results with the condition value to obtain the final
result.

\noindent\begin{minipage}[t]{0.48\textwidth}
    \begin{algorithm}[H]
        \caption{Regular conditional branching.}
    \begin{algorithmic}[1]
        \If{$\text{condition}() == 1$}
            \State $r \gets \text{foo}()$
        \Else
            \State $r \gets \text{bar}()$
        \EndIf
    \end{algorithmic}
    \label{alg:branching_classic}
\end{algorithm}
\vspace{1mm}
\end{minipage}\hfill
\begin{minipage}[t]{0.48\textwidth}
    \begin{algorithm}[H]
        \caption{Conditional branching with FHE.}
    \begin{algorithmic}[1]
        \State $c \gets \text{condition()}$
        \State $r_\text{foo} \gets \text{foo}()$
        \State $r_\text{bar} \gets \text{bar}()$
        \State $r \gets c * r_\text{foo} + (1 - c) * r_\text{bar}$
    \end{algorithmic}
    \label{alg:branching_fhe}
    \end{algorithm}
\end{minipage}

Such an approach is feasible if a small, predetermined number of branches needs to be
considered, e.g., when computing the well-known HLLC flux for the compressible Euler
equations, which has different branches based on three different wave speeds \cite{Toro09}.
It becomes intractable, however, if the number of branches is large and possibly unknown
in advance, e.g., when terminating a time step loop 
while using a dynamic, CFL-based time step size
or when considering iterative schemes with error-based stopping criteria.
Furthermore, numerical schemes that strictly require conditional branching based on current
solution values are not directly compatible with FHE. This includes methods that dynamically
modify their data structures during the simulation, such as periodic remeshing when dealing
with deforming boundaries or adaptive mesh refinement.
As a result, this poses a concrete limitation on the applicability of FHE to certain classes
of numerical simulations.

Therefore, when choosing a numerical method
for secure numerical simulations, it is important to avoid the need for control flow that
depends on encrypted data and that does not have a fixed number of branches. While the
statement that FHE can be used to compute any function (see Sec.~\ref{sec:fhe}) is not
violated, it means that in practice, the translation to an FHE algorithm is not always
feasible. As a workaround, one can consider allowing certain control flow operations to be
unencrypted, like we did in Sec.~\ref{sec:full_algorithm} for the time loop in Alg.~\ref{alg:num_sim}.

\subsection{General strategies for improving performance}
As noted in the previous subsections, most methods for dealing with
mathematical models and numerical schemes that require more complex FHE operations
incur a significant computational overhead. To mitigate this, several strategies can be
employed to optimize the performance of FHE algorithms:
\begin{itemize}
    \item \emph{Prefer parallelizable algorithms over sequential ones.} As discussed in
    Sec.~\ref{sec:multiplication}, using a binary approach for multiple multiplications can
    reduce the overall multiplicative depth of the FHE operations. This naturally applies
    also to all other functions that consume multiplicative levels. For example, consider
    loop unrolling to manually encode its iterations more efficiently.
    \item \emph{Minimize or avoid FHE-hard operations.} For instance, another
    set of state variables, e.g., using a different non-dimensionalization strategy, may save
    some (costly) division operation in each time step.
    \item \emph{Reduce the number of rotations.} While they do not consume multiplicative
    levels, they require similar computational resources as multiplications (see
    Table~\ref{tab:operations}). For example, summation over a row/column of a matrix may be
    expressed as a series of additions instead.
    \item \emph{Consider performing costly computations a priori.} For example, it might be
    more efficient to use a fixed, smaller time step based on an conservative estimate,
    rather than computing a potentially larger step size dynamically in each time step loop iteration.
\end{itemize}
A more detailed description of these and other optimization strategies can be found in
\cite{BlattGusevEtAl20}.

\subsection{Unsupported simulation types}\label{sec:unsupported_simulations}
The presented prototypes demonstrate the viability of FHE-secured numerical simulations
for some simple schemes, and the previous sections outline challenges and
solutions for extending them to more complex use cases. Some types of numerical simulations,
however, remain infeasible---either due to inherent limitations of FHE/CKKS or due to
current technical constraints.

Fundamentally unsupported are numerical methods that rely on data-dependent control flow,
such as adaptivity (adaptive time stepping, adaptive mesh refinement), nonlinear solvers
with conditional logic, event-driven simulations, or iterative schemes with
convergence-based termination. Similarly,
algorithms requiring irregular data structures or random memory access patterns are
difficult or impossible to realize with the CKKS scheme. Examples include irregular sparse
linear algebra, mesh-free particle methods, or multilevel algorithms with non-uniform
transfer operators. Finally, numerical methods that require high numerical accuracy and are
sensitive to rounding errors (typically double precision or higher), such as time
integration
for stiff PDEs, spectral methods, or eigenvalue solvers for ill-conditioned systems, are not
well-suited for use with CKKS if they require bootstrapping.

While some of these methods may be rewritten to avoid the limitations imposed by FHE through
careful algorithm design, this usually introduces an even higher computational overhead.
Thus, they represent significant challenges for many common numerical methods.

\section{Conclusions and outlook}\label{sec:conclusions}

In this work, we present the first prototype of secure numerical simulations using fully
homomorphic encryption.
We gave an overview of the CKKS scheme and its efficient
implementation within the OpenFHE library, and introduced our OpenFHE.jl and
SecureArithmetic.jl packages, which provide a user-friendly interface to OpenFHE in the
Julia programming language. A detailed analysis of the accuracy and performance of the
individual CKKS routines revealed that while each operation introduces some numerical
error, bootstrapping is by far the least accurate procedure. Similar observations were
made for the runtime behavior, where the bootstrapping mechanism is orders of magnitude
slower than other CKKS operations.

Next, we implemented the first-order upwind and second-order Lax-Wendroff
schemes for the linear advection equations for secure numerical simulations using the
SecureArithmetic.jl package. For this purpose, we introduced a $\circshift$ operation
that is compatible with the CKKS scheme and that works for both 1D and 2D data
representations. Convergence tests verified the
correctness of our implementations and showed that the CKKS encryption
does not break the convergence properties of the finite difference schemes.

We then conducted homomorphically encrypted numerical simulations by building upon the
tools and methods developed in the previous sections. Our results illustrate that it is
feasible to implement secure numerical simulations with the CKKS scheme, with results that
match the unencrypted simulations well. However, the encrypted operations introduce a large
computational overhead and are a significant source of error, to the point that we were able
to achieve less than single-precision floating point accuracy at the cost of a six-orders-of-magnitude increase in runtime. Nearly all of these undesirable effects can be attributed to
the bootstrapping procedure, which is the most inaccurate and slowest operation in the CKKS
scheme.
However, bootstrapping is necessary to refresh the multiplicative depth of the
ciphertexts, making it indispensable for long-running simulations.
Moreover, iterative bootstrapping further improves the accuracy of an FHE
algorithm, though at a non-negligible runtime penalty that needs to be taken into account.
We also showed that the CKKS scheme implementation in OpenFHE is parallelizable, and that
using multiple threads can lead to a significant reduction in the solution time.

Finally, we explored challenges---and how to address them---when extending secure
numerical simulations beyond the presented examples. Many of the discussed issues can be
mitigated through careful algorithm design and by accepting some computational overhead.
However, the lack of dynamic, value-based conditional branching precludes the
use of schemes that strictly depend on it, such as adaptive mesh refinement.

To summarize, our work demonstrates the feasibility of secure numerical
simulations using fully homomorphic encryption. We developed a user-friendly interface
to the OpenFHE library in Julia, introduced building blocks such as the
$\circshift$ algorithms, and verified that the overall simulation results agree well with
unencrypted baselines.
These contributions lay the foundation for further research in
FHE-secured scientific computing.

Looking ahead, several critical research challenges must be overcome to make secure
simulations practical at scale.
First, the numerical accuracy achievable with CKKS is inherently limited. Higher precision, e.g., through iterative bootstrapping, demands significantly more
computational resources, and double-precision accuracy appears out of reach given the
current bootstrapping implementation.
Second, the homomorphic paradigm precludes core algorithmic structures in scientific
computing, such as data-dependent control flow (necessary for iterative or adaptive schemes)
or irregular data structures (e.g., sparse matrices).
Third, the runtime overhead remains substantial---not only compared to unencrypted
simulations but also in absolute
terms. Bridging this performance gap will require improvements in CKKS implementations, new
approaches like functional bootstrapping, and the efficient use of accelerators such as GPUs
\cite{cryptoeprint:2021/508, DBLP:journals/tdsc/YangSDZLZ24, agullodomingo2025fideslib}.
Emerging hardware architectures tailored for FHE workloads, including FPGAs
\cite{cousins2017designing, cryptoeprint:2019/1066} and FHE-specific ASICs
\cite{samardzic2021f1,RPU}, could play a central role in reducing runtimes.
Addressing these limitations is the key to enabling broader use of FHE in scientific computing.
In the meantime, FHE-secured numerical
simulations only make practical sense in specialized applications where data confidentiality is of utmost importance.

Despite these open challenges, we believe that progress can be made incrementally. In particular, we plan to extend our work in several areas. We would like to increase
the complexity of the simulations by looking at systems of equations and at
three-dimensional problems.
This also necessitates an extension of the $\circshift$ operation to support multiple
ciphertexts and tensors. Furthermore, we will investigate the discretization of nonlinear
equations, which requires us to handle additional arithmetic operations such as division and
exponentiation.
Moreover, we would like to better understand the observed behavior that the accuracy of
secure numerical simulations with bootstrapping can improve over time.
Finally, we plan to continue developing the SecureArithmetic.jl
package as a tool for experimentation and rapid prototyping with encrypted algorithms, to
open up FHE to a wider audience in the computational science community.

\section*{Data availability}\label{sec:data_availability}

All code used in this work and the data generated during this study are available in our
reproducibility repository \cite{reproKholodSchlottkeLakemper24}.

\section*{CRediT authorship contribution statement}\label{sec:credit_authorship}

\textbf{Arseniy Kholod}:
Data curation, Investigation, Methodology, Software, Validation, Visualization, Writing –-
original draft, Writing --- review \& editing.
\textbf{Yuriy Polyakov}:
Methodology, Validation, Writing --- review \& editing.
\textbf{Michael Schlottke-Lakemper}:
Conceptualization, Investigation, Methodology, Resources, Software, Supervision, Validation,
Writing –- original draft, Writing --- review \& editing.

\section*{Declaration of competing interest}\label{sec:competing_interests}

The authors declare that they have no known competing financial interests or personal
relationships that could have appeared to influence the work reported in this paper.

\section*{Declaration of generative AI and AI-assisted technologies in the writing
process}\label{sec:declaration-ai}

During the preparation of this work the authors used ChatGPT in order to enhance
readability. They further used the GitHub Copilot plugin for Visual Studio Code for
assistance with mathematical notation and formatting in LaTeX. After using these
tools/services, the authors reviewed and edited the content as needed and take full
responsibility for the content of the published article.

\section*{Acknowledgments}\label{sec:acknowledgments}
We gratefully acknowledge the computing time provided by Numerical Simulation
research group of Gregor Gassner, University of Cologne.
Michael Schlottke-Lakemper further acknowledges funding through the DFG research unit FOR~5409
"Structure-Preserving Numerical Methods for Bulk- and Interface Coupling of Heterogeneous
Models (SNuBIC)" (project number 463312734), as well as through a
DFG individual grant (project number 528753982).
Moreover, we would like to thank Andrew Hazel for his valuable feedback on the
manuscript.

\appendix
\section{OpenFHE implementation details}\label{app:openfhe_details}
In the OpenFHE library, the value returned by the \texttt{GetLevel()}
function does not represent the ciphertext level as it is commonly used in literature
and in this paper. Instead, it starts at zero for a fresh ciphertext and increases by
one with each multiplication. Also, when using the OpenFHE library with the
\texttt{FLEXIBLEAUTO} scaling technique as we are here (see
Table~\ref{tab:configuration} above), the rescaling step after a ciphertext
multiplication is not performed immediately but only just before the next multiplication. Therefore, the ciphertext level reported
by \texttt{GetLevel()} function may be off by one. These implementation-specific peculiarities must be taken
into account when implementing algorithms with OpenFHE.

\bibliography{references}

\end{document}